\documentclass[a4paper,11pt]{article}
\usepackage{amssymb}
\usepackage{graphicx}
\usepackage{subfigure}
\usepackage{float}

\newcommand{\I}{{\cal I}}

\newcommand{\R}{\mathbb{R}}

\newcommand{\N}{\mathbb{N}}
\newcommand{\cH}{\mathcal{H}}

\newcommand{\cuad}{{\sqcap\kern-.68em\sqcup}}

\newcommand{\norm}[1]{\|#1\|}

\newtheorem{definition}{Definition}[section]
\newtheorem{theorem}{Theorem}[section]
\newtheorem{proposition}{Proposition}[section]

\newtheorem{lemma}{Lemma}[section]
\newtheorem{corollary}{Corollary}[section]
\newtheorem{remark}{Remark}[section]
\newcommand{\bremark}{\begin{remark} \em}
\newcommand{\eremark}{\end{remark} }

\usepackage{cite}
\begin{document}

\begin{center}{\bf  \large   On semilinear elliptic equation  with negative exponent\\[2mm]   
 arising from a closed MEMS model

 }\bigskip

 \bigskip

 {   Huyuan Chen\footnote{chenhuyuan@yeah.net}\qquad Ying Wang\footnote{yingwang00@126.com} }

\medskip
{\small  Department of Mathematics, Jiangxi Normal University,\\
Nanchang, Jiangxi 330022, PR China}  \\[10pt]
 
 {   Feng Zhou\footnote{fzhou@math.ecnu.edu.cn} }

\medskip
{\small Center for PDEs, School of Mathematical Sciences, East China Normal University,\\
Shanghai Key Laboratory of PMMP, Shanghai 200062, PR China } \\[20pt]

\begin{abstract}

This paper is concerned with the  elliptic equation  
 $-\Delta u=\frac{\lambda }{(a-u)^p}$ in a connected, bounded $C^2$  domain $\Omega$ of $\mathbb{R}^N$ 
 subject to zero Dirichlet boundary conditions, where $\lambda>0$, $N\geq 1$, $p>0$ and $a:\bar\Omega\to[0,1]$ vanishes at the boundary with the rate ${\rm dist}(x,\partial\Omega)^\gamma$
 for   $\gamma>0$.
 
 When $p=2$ and $N=2$, this equation models the closed Micro-Electromechanical Systems devices, where the elastic membrane sticks the curved ground plate on the boundary, but insulating on the boundary.  The function $a$  shapes the curved ground plate.  
 
Our aim in this paper is  to study qualitative properties of minimal solutions of this equation when $\lambda>0$, $p>0$
  and to show how the boundary decaying of $a$  works on   the minimal solutions and  the pull-in voltage. Particularly, we give a complete analysis for the stability of the minimal solutions.

 \end{abstract}

  \end{center}

 % \tableofcontents \vspace{1mm}
  \vspace{1mm}


\vspace{2mm}
%%%%%%%%%%%%%%%%%%%%%%%%%%%%%%%%%%%%%%%%%%%%%%%%%%%%%%%%%%%%%%%%%%%%%%%%%%%%%%%%%%%%%%%%%%%%%%%%%%%%%%%%%%%%%%%%%%%%%%%%%%
%%%%%%%%%%%%%%%%%%%%%%%%%%%%%%%%%%%%%%%%%%%%%%%%%%%%%%%%%%%%%%%%%%%%%%%%%%%%%%%%%%%%%%%%%%%%%%%%%%%%%%%%%%%%%%%%%%%%%%%%%%

\setcounter{equation}{0}
\section{Introduction}

Our purpose of this paper is to consider the  minimal solutions of
  elliptic equation
\begin{equation}\label{eq 1.1}
\left\{\arraycolsep=1pt
\begin{array}{lll}
\displaystyle -\Delta    u =  \frac{\lambda}{(a-u)^{p}}\quad  &{\rm in}\quad\ \Omega,
 \\[2.5mm]
 \phantom{- }
\displaystyle  0<u<a\quad &{\rm in}\quad\ \Omega,
 \\[1mm]
 \phantom{-\Delta   }
\displaystyle   u=0\quad &{\rm on}\quad   \partial \Omega,
\end{array}
\right.
\end{equation}
where $\lambda,\, p>0$, $\Omega$ is a connected, bounded $C^2$ domain in $\R^N$ with $N\geq 1$ and $a:\bar \Omega\to [0,1]$  is a  function in $  C_0^2 (\Omega) :=C_0 (\Omega)\cap C^2(\Omega)  $
being positive inside and vanishing  on the boundary.

When  $p=2$ and $a\equiv1$,   problem (\ref{eq 1.1}) models the standard   Micro-Electromechanical Systems (MEMS), which    are   used to combine electronics with micro-size mechanical devices in the design of various types of microscopic machinery.  The key component in  MEMS devices is called the electrostatic actuation, which
is based on an electrostatic-controlled tunable, it is a simple idealized electrostatic device. The upper part of this electrostatic device consists of a thin and deformable elastic membrane that is fixed along its
boundary and which lies above a rigid grounded plate. This elastic membrane is modeled as a dielectric with
a small but finite thickness. The upper surface of the membrane is coated with a negligibly thin metallic
conducting film.   The unknown function  $u$ in (\ref{eq 1.1}) stands for the deformation of the elastic membrane and  the term on the right hand side of the equation  is the Coulomb force.

When a voltage $\lambda$ is applied to the conducting film, the thin dielectric membrane deflects
towards the bottom plate, and when $\lambda$ is increased beyond a certain critical value $\lambda^*$, known as pull-in
voltage, the steady state of the elastic membrane is lost. This  proceeds to touchdown or snap through at a finite time creating the so-called pull-in instability.  Equations related to MEMS  have been  attracting great attentions in the last decades. 
  The authors  \cite{EGG,GG,G,GPW,GZZ} studied MEMS equation 
\begin{equation}\label{eq 1.1241}
 -\Delta    u =  \frac{\lambda f(x)}{  (1-u)^{2}} \ \ {\rm in}\ \, \Omega,\qquad u=0\ \ {\rm on}\ \, \partial\Omega,
\end{equation}
  which  models a simple electrostatic
MEMS device consisting of a thin dielectric elastic membrane with boundary supported at $0$ above a
rigid ground plate located at $1$.  Here the function $f$  represents the permittivity profile.  
Qualitative properties of these type equations  have been studied in \cite{ES,ES1,Li,LW,LW1,LY,P,YZ} and the references therein. When $f$ is replaced by $1+|\nabla u|^2$,  the   MEMS equations with   fringing fields are discussed in \cite{DW,GH}.

  Recently, the authors in \cite{CWZ} used the equation
\begin{equation}\label{eq 1.1-2}
 -\Delta    u =  \lambda  (a-u)^{-2}\ \  {\rm in}\ \, \Omega,
 \qquad u=0\ \  {\rm on}\ \,   \partial \Omega
\end{equation}
to model  a closed MEMS device
that the static deformation  of the surface of membrane when it is applied  voltage $\lambda$,
where $a$ is the shape of the ground plate, which is no longer plat. The elastic membrane contacts the ground plate along the boundary, but insulating on the boundary, see the Picture 1 and Picture 2 below.    Let   $a:\bar\Omega\to[0,1]$ be in the class of $C^2_0(\Omega)$ and satisfy that 
\begin{equation}\label{1.1}
\frac1{\kappa}\rho(x)^\gamma \leq a(x) \leq \kappa \rho(x)^\gamma,\quad \forall \, x\in\Omega
\end{equation}
for some $\kappa\geq 1$, $\gamma>0$ and 
$$\rho(x)=\min\Big\{\frac12, {\rm dist}(x,\partial\Omega)\Big\}.$$ 
  \begin{figure}[H]
   \centering
   \subfigure{
   \begin{minipage}{65mm}
   \includegraphics[scale=0.55]{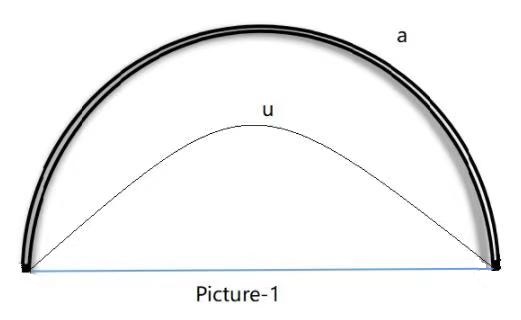}
  \end{minipage}
  \begin{minipage}{65mm}
   \includegraphics[scale=0.5]{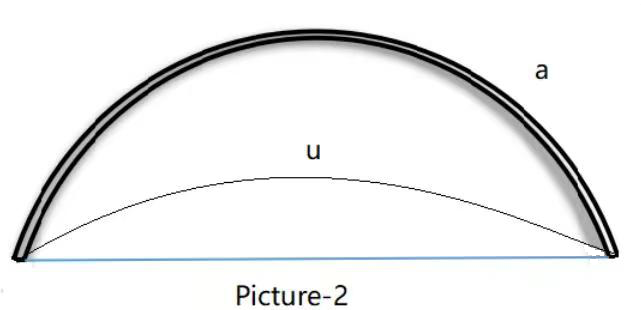}
  \end{minipage}
  }

   {\footnotesize  The shape of ground plates have the behaviors $a\sim \rho^{\frac12}$ in Picture 1
   and $a\sim\rho$ in Picture 2.  }
  \end{figure} 
\noindent It is proved in  \cite{CWZ} that for   $\gamma\in(0,\frac23]$  there exists a critical value
$\lambda^*$  (pull-in voltage) depending on $\kappa, \gamma$ such that if $\lambda\in(0,\lambda^*)$,  problem (\ref{eq 1.1-2}) has a minimal solution, while for $\lambda>\lambda^*$, no solution  exists for (\ref{eq 1.1-2}).

 Note that the nonlinearity with negative quadratic stands  the field force. 
While   different types of forces would be presented by different nonlinearities, such as  
$ (1-u)^{-2}$ stands for the Coulomb force, negative quartic $ (1-u)^{-4}$ stands for  effects of Casimir force, see the references \cite{BPS,L,La},  MEMS with the fringing field \cite{DW} and with general negative exponent refers to \cite{ZL}. 
  Our purpose of this paper is to consider the qualitative property of the minimal solutions of the closed MEMS equation (\ref{eq 1.1}) with the  nonlinearity  having  a negative exponent.

 For this, we have to derive the critical exponent depending on the decaying rate $\rho^{\gamma}$ of $a$ near the boundary and then to find out pull-in
voltage of $\lambda$  such that   (\ref{eq 1.1}) admits  a minimal solution.    For $\gamma\in(0,1]$, we denote 
 \begin{equation}\label{crit 1}
 p^*_\gamma= \frac{2}{\gamma}-1,
 \end{equation}
 and
 \begin{equation}\label{inter I}
\I_\gamma=\left\{\arraycolsep=1pt
\begin{array}{lll}
\displaystyle  (0, p^*_\gamma] \quad  &{\rm if}\ \   \gamma\in(0,1),
 \\[2mm]
 \phantom{ }
\displaystyle (0,\, 1) \quad &{\rm if}\ \   \gamma=1.
\end{array}
\right. 
\end{equation}

  Our main result on the existence  states as following.

\begin{theorem}\label{teo 1}
Assume that    $a\in C_0^2(\Omega)$ satisfies (\ref{1.1}) with 
$\kappa\geq 1$, $\gamma\in(0,1]$ and $p\in \I_\gamma$. 

Then   there exists a  pull-in voltage $\lambda^*_{\gamma,p}\in(0,+\infty)$ depending on $\gamma,\, p$ such that\smallskip 

$(i)$ for $\lambda\in(0,\lambda^*_{\gamma,p})$, problem (\ref{eq 1.1}) admits a minimal solution $u_{p,\lambda}$ and the mapping: $\lambda\mapsto u_{p,\lambda}$
is   strictly  increasing and continuous.

Furthermore,  for some  $c_1\ge 1$.
\begin{equation}\label{e 1.1-sub}
 \arraycolsep=1pt
\begin{array}{lll} 
  \frac{\lambda}{c_1}\rho(x)^{\min\{1,2-p\gamma\}}\le u_{p,\lambda}(x)\le c_1 \lambda \rho(x)^{ \gamma },\quad \forall\, x\in
\Omega,\quad  {\rm if}\quad p\not=\frac1\gamma\\[3mm] 
  \frac{\lambda}{c_1}\rho(x) \ln \frac1{\rho(x)}  \le u_{p,\lambda}(x)\le c_1 \lambda \rho(x) ^{ \gamma }  ,\quad\ \forall\, x\in \Omega,\qquad \    {\rm if}\quad p=\frac1\gamma.
\end{array}
\end{equation}

$(ii)$ for $\lambda>\lambda^*_{\gamma,p}$, there is no solution for (\ref{eq 1.1}).\smallskip\smallskip

Moreover, there hold \\[1mm]
$(iii)$ for $\gamma\in(0,1]$ fixed,  the mapping $p\mapsto \lambda^*_{\gamma,p}$ is decreasing and 
$$ \lim_{p\to0^+}\lambda^*_{\gamma,p}\leq \frac{\|a\|_{L^1(\Omega)}}{\|\mathbb{G}_\Omega[1]\|_{L^1(\Omega)}}, $$ 
where $\mathbb{G}_\Omega[1]$ is the solution of 
$$-\Delta u=1\quad{\rm in}\ \,  \Omega,\quad\ \  u=0\quad{\rm on}\ \, \partial\Omega; $$
$(iv)$ if $\gamma=1$, $\lambda^*_{\gamma,p}\to0$ as $p\to 1^-$.
\end{theorem}

 In the super critical case, we have the following nonexistence results.

\begin{theorem}\label{teo 2}
Assume that $a\in  C_0^2 (\Omega) $  satisfies that  for $\gamma>0$
\begin{equation}\label{1.1-sup}
0<a(x) \leq   \kappa \rho(x)^\gamma,\quad \forall \, x\in\Omega.
\end{equation}
 Then  problem (\ref{eq 1.1}) admits no nonnegative solution for any $\lambda>0$ if one of the following holds: \smallskip

  $(i)$ \ \ $p>p^*_\gamma$ for $\gamma\in(0,1)$; \smallskip
  
$(ii)$ \ $p\geq 1$  \, for $\gamma=1$;\smallskip

$(iii)$ $p>0$ \, for $\gamma>1$.

\end{theorem}

We obtain optimal ranges: $p\in(0,p^*_\gamma]$ if $\gamma\in(0,1)$ and  $p\in(0,1)$ if $\gamma=1$ for the existence of  pull-in voltage $\lambda^*_{\gamma,p}$ and the minimal solution $u_{p,\lambda}$ of (\ref{eq 1.1}) for $\lambda\in(0,\lambda^*_{\gamma,p})$  from Theorem \ref{teo 1} and the pull-in voltage $\lambda^*_{\gamma,p}=0$ under  the setting of Theorem \ref{teo 2}. 
 \begin{figure}[H]
   \centering
   \subfigure{
   \begin{minipage}{65mm}
   \includegraphics[scale=0.35]{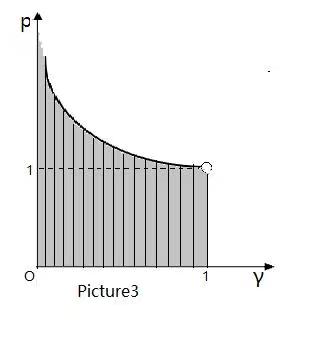}
  \end{minipage}
  }

   {\footnotesize  The shade region of coordinate $(\gamma,\, p)$ represents the existence of classical minimal solutions of (\ref{eq 1.1}),  the above curve of the shade region, except a point $(1,1)$ does the existence of the extremal solution. }
  \end{figure}

Normally, it is a challenging task to  consider the limit of the minimal solutions $\{u_{p,\lambda}\}_{\lambda}$ as $\lambda\to\lambda^*_{\gamma,p}$, which is called  as the extremal solution.  Indeed, the extremal solution  used to be found in the distributional sense and then to improve  the regularity until
  to the classical sense when the dimension $N$ co-works with parameters $\gamma$ and $p$  suitably.   
  
\begin{definition}\label{def 1}
A function $u$ is a weak solution of (\ref{eq 1.1}) if $0\le u\le a$ and $(a-u)^{-p}\in L^1_{loc}(\Omega)$
$$
\int_\Omega u (-\Delta) \xi dx=\int_\Omega  \lambda(a-u)^{-p} \xi  dx,\quad \forall\, \xi\in C_c^2(\Omega),
$$
where  $C_c^2(\Omega)$ is the set of all function  in $C^2(\Omega)$ with compact support in $\Omega$.

A solution (or weak solution) $u$ of (\ref{eq 1.1}) is stable (resp. semi-stable) if %$(a-u)^{-1-p}\in L^1(\Omega,\rho^2 dx)$
$$
\int_\Omega \frac{p\lambda \varphi^2}{(a-u)^{p+1}} dx<\int_\Omega  |\nabla\varphi|^2 dx,\quad ({\rm resp.}\ \le)\quad \forall\, \varphi\in \cH^1_0(\Omega)\setminus\{0\},
$$
where $\cH_0^1(\Omega)$ is the Hilbert space with the inner product 
$\displaystyle \langle w,v\rangle= \int_\Omega \nabla w\cdot \nabla v dx$ and the induced norm $\norm{v}=\sqrt{\langle v,v\rangle }$. 
\end{definition}

In the extremal case, we have the following existence results.

\begin{theorem}\label{teo 5}
Assume that $a\in   C^2_0(\Omega)$ satisfies (\ref{1.1}) with 
$\kappa\geq 1$ and $\gamma\in(0,1]$. Let   $\lambda^*_{\gamma,p}$ be given in Theorem \ref{teo 1} for $p\in \I_\gamma$.  Then  problem (\ref{eq 1.1}) with $\lambda=\lambda^*_{\gamma,p}$ admits a weak solution  $ u_{p,\lambda^*_{\gamma,p}}$.  Moreover, for any $\beta\in(0,\gamma)$,  there exists $c_\beta>0$ such that
$$% \begin{equation}\label{4.0.5}
 \norm{u_{p,\lambda^*_{\gamma,p}}}_{W^{1,\frac{N}{N-\beta}}(\Omega)}\le c_\beta
$$% \end{equation}
and
$$
 \int_{\Omega} \frac{\rho^{1-\beta}}{(a-u_{p,\lambda^*_{\gamma,p}})^p} dx\le c_\beta.
$$

\end{theorem}

Our  aim is to  consider the stability of minimal solutions and the regularity of extremal solutions.

\begin{theorem}\label{teo 3}
Assume that $a\in   C^2_0(\Omega)$ satisfies (\ref{1.1}) with 
$\kappa\geq 1$ and $\gamma\in(0,1]$. For $p\in \I_\gamma$ and $\lambda\in (0,\lambda^*_{\gamma,p}]$, let  $u_{p,\lambda}$ be the minimal solution of (\ref{eq 1.1}). 
 Then    for  $\lambda\in (0,\lambda^*_{\gamma,p})$
the minimal solution  $u_{p,\lambda}$ is  stable 
  and  the minimal solution $ u_{p, \lambda^*_{\gamma,p}}$ is  semi-stable.

\end{theorem}
 
 Finally,  we try to improve the  regularity of extremal solution via the stability.
 To this end, let us introduce some critical exponents $p^\#_{_N}$ and $ q^\#_{_N}$.
Denote   
 \begin{equation}\label{ex p}
p^\#_{_N}=\left\{\arraycolsep=1pt
\begin{array}{lll}
\displaystyle \frac{1}{(\sqrt{N-2}-1)^2-1} \quad  &{\rm if}\ \  N\geq7,
 \\[3mm]
 \phantom{ }
\displaystyle +\infty \quad &{\rm if}\ \   N\in\{1,2,3,4,5,6\}
\end{array}
\right. 
\end{equation}
and
 \begin{equation}\label{ex q}
q^\#_{_N}=\left\{\arraycolsep=1pt
\begin{array}{lll}
\displaystyle \frac{1}{(\sqrt{\frac N2-2}-1)^2-1} \quad  &{\rm if}\ \  N\geq13,
 \\[3mm]
 \phantom{ }
\displaystyle +\infty \quad &{\rm if}\ \   N\in\{1,2,\cdots,12\}.
\end{array}
\right. 
\end{equation}
 When $N\geq7$, $p^\#_{_N}$ is the unique zero point in $[0,+\infty)$ of  $f_0(t)-N=0$, where
 \begin{equation}\label{fff}
 f_0(t):=3+\frac1t  + 2\sqrt{1 +\frac1t},\quad \forall\, t>0, 
 \end{equation}
 which is strictly decreasing and $\displaystyle \lim_{t\to+\infty} f_0(t)=5$. 
 When $N\geq13$,  $q^\#_{_N}$ is the unique zero point in $[0,+\infty)$ of  $f_0(t)-\frac N2=0$.
 If $p^\#_{_N}<p^*_\gamma$, then for $p^\#_{_N}\leq p< \min\{p^*_\gamma, q^\#_{_N}\}$, we have that 
 $$\frac N2<f_0(p)\leq N.$$ 
 
\begin{theorem}\label{teo 6}
Assume that $a\in   C^2_0(\Omega)$ satisfies (\ref{1.1}) with 
$\kappa\geq 1$ and $\gamma\in(0,1]$. For $p\in(0,\, p^*_\gamma]$, let  $u_{p,\lambda}$ be the minimal solution of (\ref{eq 1.1}) with $\lambda\in (0,\lambda^*_{\gamma,p})$.\smallskip

Then $(i)$  %  Assume more that
%\begin{equation}\label{4.0.6-1}
%\int_\Omega |\Delta a|\, dx <+\infty.
% \end{equation}
  for $$p\in \I_\gamma\quad{\rm and}\quad  p<p^\#_{_N},$$  
the extremal solution $u_{p,\lambda^*_{\gamma,p}}$ is a classical solution of (\ref{eq 1.1}).\smallskip

$(ii)$  Let %(\ref{4.0.6-1}) hold true and 
$p^\#_{_N}<p^*_\gamma$, for   
\begin{equation}\label{4.0.6-2}
 p^\#_{_N}\leq p< \min\{p^*_\gamma, q^\#_{_N}\}  \quad  {\rm if}\ \  N\geq 7,
 \end{equation}
then
$ u_{p,\lambda^*_{\gamma,p}} \in C^\alpha_{loc}(\Omega)\cap C_0(\Omega)$
for 
$$\alpha=2-\frac{N}{f_0(p)}\in(0,1],$$
where and in the sequel we use the notation $C^{\sigma}=C^{0,1}$ if $\sigma=1$, $C^{\sigma}=C^{1,\sigma-1}$ if $\sigma\in(1,2)$. 
\end{theorem}

%We then have that $p^\#_\gamma<\frac13<1\leq p^*_\gamma$ and 
%$f(p^*_\gamma)=\frac{2(2-\gamma)^2(3-\gamma+\sqrt{2(2-\gamma)})}{\gamma^2}.$
%Note that  $f(p^*_1)=f(1)=2(2+\sqrt{2})$ and $f(p^*_\gamma)\to+\infty$ as $\gamma\to0^+$.

%Combine the following properties of $p^*_\gamma, \lambda^*_{\gamma,p}$ from following Proposition, we use Picture 2 and Picture 3 to illustrate the regions of $(p,\lambda)$ for the stabilities of minimal solutions of (\ref{eq 1.1}).

Due to the decay of function $a$, the study of extremal solutions becomes   subtle.  The difficulty  to obtain the stability of $u_{p,\lambda}$ comes from the blowing up of $\frac{1}{(a-u_{p,\lambda})^{p+1}}$ at the boundary.
\begin{remark}
$(i)$ When $p=2$,  \cite{CWZ}  obtains the stability for only $\lambda\in(0,\lambda_*)$, where $\lambda_*\leq \lambda^*_{\gamma,p}$ and only getting the equal  for  the particular case $\gamma=\frac23$.
 From  Theorem \ref{teo 3}  we give a complete stability for all minimal solutions in a general model.
 
 $(ii)$ Involving the critical  exponents $p^\#_{_N}$ and $q^\#_{_N}$, our model could give more information for the regularity of the extremal solutions.  For $p<p^\#_{_N}$, thanks to the regularity of  extremal solutions, it provides to the possibility to construct a second solution for $\lambda<\lambda^*_{\gamma,p}$. For $a\equiv 1$ and $p=2$,  multiple solutions are
 constructed in \cite{EGG,GG}.
 
 $(iii)$ Thanks to the boundary decay, our stability of minimal solutions in Theorem  \ref{teo 5}  are based  on the convergence as $\epsilon\to0^+$,  of minimal solutions  of 
 $$\left\{\arraycolsep=1pt
\begin{array}{lll}
\displaystyle -\Delta    u =  \frac{\lambda}{(a+\epsilon-u)^{p}}\quad  &{\rm in}\quad\ \Omega,
 \\[3mm]
 \phantom{- }
\displaystyle  0<u<a+\epsilon\quad &{\rm in}\quad\ \Omega,
 \\[2mm]
 \phantom{-\Delta   }
\displaystyle   u=0\quad &{\rm on}\quad   \partial \Omega
\end{array}
\right.
$$
for $\epsilon\in(0,1)$. %This convergence is also  essential for the derivation of our stability of minimal solutions and extremal solutions.  

 $(iv)$ The critical exponent $p^\#_{_N}$ could be seen explicitly 
 \begin{center}
\renewcommand{\arraystretch}{1.6}
    \begin{tabular}{ l | l | l | l | l | l | l }
   \hline $N$ & 6 & 7 & 8 & 9&10&11 \\ \hline
    $p^\#_{_N}$ & $+\infty$& $1+\frac25\sqrt{5}\in(2,3)$&$\frac12+\frac{\sqrt{6}}{6}\in(\frac12,1)$&$\frac13+\frac{2\sqrt{7}}{21}\in(\frac12,1)$&$\frac14+\frac{2\sqrt{2}}{8}\in(\frac13,\frac12) $&$\frac13$
\\[1mm]  \hline
\end{tabular}
\renewcommand{\arraystretch}{1}
\end{center}
Note that  for $p=2$ the nonlinearity represents the Coulomb force,  the critical exponent  $p^\#_{_N}> 2$  if  $N\leq 7$ and $p^\#_{_N}< 2$ for $N\geq 8$. 
 \end{remark}
 
 Also we'd like to point out that if $u_{p,\lambda^*_{\gamma,p}}$ is  regular or $C^2$ in $O\Subset \Omega$,  then it does not contact $a$ in $O$.  Inversely, thanks to the singularity of $(a-u)^{p}$ at the contracting points, the regularity of the extremal solution represents how the elastic membrane contacts the ground plate. 
   {\it It is still open what the regularity of the extremal solutions is for $p>q^\#_{_N}$. }

 \smallskip

The rest of the paper is organized as follows.  Section 2 is devoted to   the existence of the pull-in voltage $\lambda^*_{\gamma,p}$ such that problem (\ref{eq 1.1}) admits a minimal solution when $\lambda\in(0,\lambda^*_{\gamma,p})$ and 
the existence of weak minimal extremal solutions.  In Section  3, we do show the stability for the minimal solutions and semi-stability of extremal solutions. Finally
 we  improve  of the regularity for  extremal solutions in Section 4.

\setcounter{equation}{0}
\section{ Existence }

\subsection{  Minimal solutions }
Denote by $G_\Omega$ the Green kernel of $-\Delta$ in $\Omega\times\Omega$ under the zero boundary condition and by $\mathbb{G}_\Omega[\cdot]$ the
Green operator defined as
$$\mathbb{G}_\Omega[f](x)=\int_{\Omega} G_\Omega(x,y)f(y)dy,\quad \forall\, f\in L^1(\Omega,\rho dx). $$
To derive the existence of  a minimal solution of problem (\ref{eq 1.1}), the following estimates play an important role.
\begin{lemma}\label{lm 2.1}\cite[Lemma 2.1]{CWZ}
Let $\tau\in(0,2)$, $A_{\frac12}=\big\{x\in\Omega: \ \rho(x)<\frac14\big\}$ and  for $x\in A_{\frac14}$ 
\begin{equation}\label{varrho}
\varrho_\tau(x)=\left\{ \arraycolsep=1pt
\begin{array}{lll}
 \rho(x)^{\min\{1,\tau\}}\ \ &{\rm if}\ \ \tau\in (0,1)\cup(1,2),
 \\[2mm]
\rho(x)  \ln \frac1{\rho(x)} \ & {\rm if}\ \ \tau=1,
\end{array}
\right.
\end{equation}
where we recall that $ \rho(x)=\min\{\frac12, {\rm dist}(x,\partial\Omega)\}$.

Then $\varrho_\tau\in C_{loc}^{0,1}$ and there exists $c_\tau>1$ such that
$$\frac1c_\tau\varrho_\tau(x) \le \mathbb{G}_\Omega[\rho^{\tau-2}](x)\le c_\tau\varrho_\tau(x),\quad\forall x\in  \Omega.$$

\end{lemma}

 Now we are ready to show the existence of pull-in voltage $\lambda^*_{\gamma,p}$ to problem (\ref{eq 1.1}) such that
 (\ref{eq 1.1}) admits a solution for $\lambda\in(0,\lambda^*_{\gamma,p})$ and no solution for $\lambda>\lambda^*_{\gamma,p}$. \medskip

\noindent {\bf Proof of Theorem \ref{teo 1}.}  {\it Existence of minimal solution. }  Let $v_0\equiv0$ in $\bar\Omega$ and
$$v_1=\lambda \mathbb{G}_\Omega[ a^{-p}]>0,$$
then by (\ref{1.1}) and Lemma \ref{lm 2.1}, we have that
\begin{equation}\label{2.428}
 v_1  = \lambda \mathbb{G}_\Omega[ a^{-p}] \le  \frac{\lambda}{\kappa^p} \mathbb{G}_\Omega[\rho^{-p\gamma}]  \le  \frac{ \lambda}{\kappa^p}c_{2} \varrho_{ _{2-2\gamma}},
\end{equation}
where $c_{2}>0$ depending on $\gamma$ and $\varrho_{_{2-2\gamma} }$ is given by (\ref{varrho}).

From the natural constraint $0<u<a$ in $\Omega$ in (\ref{eq 1.1}), we need  
$$2-p\gamma\geq \gamma,$$
which holds for $p\in(0,p^*_\gamma]$.  If $p\not=\frac1\gamma$, we observe that $\min\{1,2-p\gamma\}\ge \gamma$ and by (\ref{2.428}),
$$
v_1(x)\le  c_{3}\frac{ \lambda}{\kappa^p} \rho(x)^{\min\{1,2-p\gamma\}}\le c_{3}\frac{\lambda}{\kappa^p}\rho(x)^\gamma,\quad\forall\, x\in\Omega.
$$
When   $p=\frac1\gamma$, we see that $2-p\gamma=1$ and  if $\gamma<1$,  (\ref{2.428}) implies that 
$$
v_1(x)\le c_{3} \frac{ \lambda}{\kappa^p} \rho(x) \ln\frac1{\rho(x)} \le \frac{\lambda}{\kappa^p}c_4\rho(x)^\gamma,\quad\forall\, x\in\Omega
$$
that is,
$$v_1(x)\le c_4\frac{\lambda}{\kappa^p}\rho(x)^\gamma,\quad\forall\, x\in\Omega,$$
where $c_4>0$ depends on $\gamma$, but independent of $\lambda$ and $\kappa$.

Fixed $\mu\in(0, \frac1\kappa)$,  let $\lambda_1=\frac{ \kappa^{p}\mu}{c_4} >0$ such that
%\begin{equation}\label{2.1}
$$c_4 \frac{\lambda}{\kappa^p} \leq \mu<\frac1\kappa, \ \ \, \forall\,\lambda\leq \lambda_1,$$
then  for any $\lambda\leq \lambda_1$,
$$v_1(x)\le \mu\rho(x)^\gamma,\quad \forall\, x\in\Omega.$$
%\end{equation}

Let $v_2=\lambda \mathbb{G}_\Omega[ (a-v_1)^{-p}]$, by the fact that $a(x)\ge \frac1\kappa\rho(x)^\gamma>\mu\rho(x)^\gamma\ge v_1(x)>0$ for $x\in\Omega$ and and Lemma \ref{lm 2.1},
 we have that  
$$v_1=\lambda\mathbb{G}_\Omega[ a^{-p}]\le v_2\le  \frac{\lambda}{(\kappa^{-1}-\mu)^p}\mathbb{G}_\Omega[\rho^{-p\gamma}]
 \le  \frac{c_4 \lambda}{(\kappa^{-1}-\mu)^p}\rho^\gamma\quad {\rm in}\ \ \Omega.$$
 Let
 \begin{equation}\label{def mu}
 \mu^*=\frac{\kappa^{-1}}{p+1} \quad\ {\rm and}\quad\   \lambda_\#=(\kappa^{-1}-\mu^*)^p\mu^*=  \frac1{c_4}  p^p \big(\frac{k^{-1}}{p+1}\big)^{p+1}  ,
 \end{equation} 
 then $\mu^*$  maximizes the function of $(\kappa^{-1}-\mu)^p\mu$ for $\mu\in(0,\, \kappa^{-1})$ and 
$$\frac{c_4\lambda}{(\kappa^{-1}-\mu^*)^p}\le \mu^*<\kappa^{-1}, \ \ \ \forall\, \lambda\leq\lambda_\#.$$
Then for $\mu=\mu^*$ and any $\lambda\leq \lambda_\#$, we obtain that
$$v_2(x)\le \mu^*\rho(x)^\gamma,\quad \forall\, x\in\Omega.$$

Iterating above process, we have that for $\lambda\in(0, \lambda_\#]$,
 $$v_n:=\lambda\mathbb{G}_\Omega[(a-v_{n-1})^{-p}]\ge v_{n-1},\quad\forall\, n\in\N$$
and
$$
v_n(x)\le \mu^*\rho(x)^\gamma,\quad\forall\, x\in\Omega.
$$
Thus, the sequence $\{v_n\}_n$ converges, denoting by 
$$u_{p,\lambda}=\lim_{n\to\infty} v_n,$$
 then $u_{p,\lambda}$ is  a classical solution of (\ref{eq 1.1}).

We claim that $u_{p,\lambda}$ is the minimal solution of  (\ref{eq 1.1}), that is, for any positive solution $u$ of (\ref{eq 1.1}), we always have
$u_{p,\lambda}\le u$. Indeed, there holds $u\ge v_0$ and then
$$u=\lambda \mathbb{G}_\Omega[ (a-u)^{-p}]\ge\lambda\mathbb{G}_\Omega[ a^{-p}]=v_1.$$
One can show inductively that $$u\ge v_n\quad{\rm for\ all\ }\ n\in\N \quad\Longrightarrow   \quad u\ge u_{p,\lambda}.$$

From the approximating sequence  $\{v_n\}_n$ and $a\leq 1$, we can derive directly that 
 the mapping $\lambda\mapsto u_{p,\lambda}$ is increasing and the mapping $p\mapsto u_{p,\lambda}$ is decreasing, thanks to $0<a-u<a\leq 1$.
 Moreover, if problem  (\ref{eq 1.1})
 has a  super solution $u$  for $\lambda=\lambda_0>0$,
then (\ref{eq 1.1}) admits a minimal solution $u_{p,\lambda}$   for all $\lambda\in(0,\lambda_0]$, since at this moment function $u$ provides
an upper bound for the increasing sequence $\{v_n\}_n$.

Let us denote 
   $$\lambda^*_{\gamma,p}=\sup\big\{\lambda>0:\,  (\ref{eq 1.1}) \ {\rm has \  a \  minimal \  solution \ with\ such} \ \lambda \big\}.$$
Obviously, $\lambda^*_{\gamma,p}\geq \lambda_\#>0$.  

 For $0<\lambda_1<\lambda_2<\lambda^*_{\gamma,p}$, we know that $0\le u_{p,\lambda_1}\le u_{p,\lambda_2}\le a$ in $\Omega$, then
\begin{eqnarray}
 -\Delta (u_{p,\lambda_2}-u_{p,\lambda_1}) &=& \frac{\lambda_2}{(a-u_{p,\lambda_2} )^p} - \frac{\lambda_1}{(a-u_{p,\lambda_1} )^p}\nonumber\\
    &\ge & \frac{\lambda_2-\lambda_1}{(a-u_{p,\lambda_1} )^p}\ge \frac{\lambda_2-\lambda_1}{a^p}>0,\label{ooo-1}
\end{eqnarray}
which implies that
 \begin{equation}\label{e 1.1}
 u_{p,\lambda_2}-u_{p,\lambda_1}\ge (\lambda_2-\lambda_1)\mathbb{G}_\Omega[a^{-p}]>0.
 \end{equation}
It infers that $u_{p,\lambda}<a$ in $\Omega$ for any $\lambda<\lambda^*_{\gamma,p}$ and
by the standard interior regularity results, we have that $u_{p,\lambda}\in C_0^{2}(\Omega)$. Moreover, we observe that the map $\lambda\in(0,\lambda^*_{\gamma,p})\mapsto u_{p,\lambda}$ is strictly increasing.

\smallskip

The continuity of   the map $\lambda\in(0,\lambda^*_{\gamma,p})\mapsto u_{p,\lambda}$ follows by Proposition \ref{pr 4.2-con}, which is addressed in Section \S 3 below. 

\smallskip

{\it Upper estimate of $\lambda^*_{\gamma,p}$ and non-existence. }  We prove that $\lambda^*_{\gamma,p}<+\infty$.  If not, then for any $\lambda>0$, problem  (\ref{eq 1.1}) has the minimal solution $u_\lambda$.
Let $A_{\delta}=\{x\in\Omega:\, \rho(x)<\delta\}$ and $n\in\N$,
\begin{equation}\label{101}
\eta_n=1 \ \ {\rm in} \  \Omega\setminus{A_{1/n}}, \quad  \   \eta_n=0 \ \ {\rm in} \  A_{1/{2n}},\quad  \  \eta_n\in C^2(\Omega),
\end{equation}
and
$\xi_n= \mathbb{G}_{\Omega}[1]\eta_n$ in $\Omega$, we observe that $\xi_n\in C_c^2(\Omega)$ and
\begin{eqnarray*}
\int_\Omega u_{p,\lambda} (-\Delta \xi_n) dx=\int_\Omega u_{p,\lambda} (\eta_n-\nabla \mathbb{G}_{\Omega}[1]\cdot\nabla\eta_n+(-\Delta \eta_n)\mathbb{G}_{\Omega}[1])dx.
\end{eqnarray*}
By the fact that $|\nabla \eta_n|\le c_5n$, $|(-\Delta)\eta_n|\le c_6n^2$ and
$\mathbb{G}_{\Omega}[1]\le c_7\rho$ in $\Omega$, then we have that
\begin{equation}\label{eq9281}
\int_\Omega u_{p,\lambda} (-\Delta\xi_n) dx\le\int_\Omega u_{p,\lambda} dx+c_5n\int_{A_{1/n}}u_{p,\lambda} dx + c_8n^2\int_{A_{1/n}}u_{p,\lambda} \rho dx
\le c_{9},
\end{equation}
where $c_5, \cdots, c_9 >0$ are independent on $n$ and $A_{t}=\{x\in \Omega: {\rm dist}(x,\partial\Omega)<t\}$.

Since $u_{p,\lambda}$ is the minimal solution of  (\ref{eq 1.1}), then
$$\int_\Omega u_{p,\lambda} (-\Delta\xi_n) dx=\int_\Omega \nabla u_{p,\lambda}\cdot \nabla{\xi_n} dx
=\int_\Omega (-\Delta u_{p,\lambda})\xi_n dx=\int_\Omega \frac{\lambda \xi_n}{(a-u_{p,\lambda})^p} dx.$$
Passing to the limit as $n\to\infty$ and combining with (\ref{eq9281}), we see that
$$\int_\Omega \frac{\lambda \mathbb{G}_{\Omega}[1]}{(a-u_{p,\lambda})^p} dx\le c_{10}.$$
Thus,
\begin{eqnarray*}
\int_\Omega a(x) dx\ge \int_\Omega u_{p,\lambda}(x) dx&=&\int_\Omega \mathbb{G}_\Omega[1]  (-\Delta u_{p,\lambda})   dx\\
 &=& \lambda \int_\Omega \frac{\mathbb{G}_\Omega[1](x)}{[a(x)-u_{p,\lambda}(x)]^p}dx
  >  \lambda \int_\Omega \mathbb{G}_\Omega[1] a^{-p} dx,
\end{eqnarray*}
which implies that
$$
\lambda<  \frac{\int_\Omega a dx }{\int_\Omega \mathbb{G}_\Omega[1] a^{-p}dx}.
$$
Therefore,
$$
\lambda^*_{\gamma,p}\leq  \frac{\int_\Omega a(x) dx }{\int_\Omega \mathbb{G}_\Omega [1] a^{-p}dx}.
$$
Here $ \mathbb{G}_\Omega [1] a^{-p}$ is integrable for $p\in(0,\frac{2}{\gamma})$, where 
$\frac{2}{\gamma}>p^*_\gamma$ and
$$\lim_{p\to0^+}\frac{\int_\Omega a  dx }{\int_\Omega \mathbb{G}_\Omega[1] a^{-p}dx}=\frac{\int_\Omega a  dx }{\int_\Omega \mathbb{G}_\Omega[1]  dx}.$$
\smallskip

Let $p_1\leq p_2<p^*_\gamma$ and then for any $\lambda\in(0,\lambda^*_{\gamma,p_2})$
$$-\Delta u_{p_2,\lambda}=(a-u_{p_2,\lambda}) ^{-p_2}\geq (a-u_{p_2,\lambda})^{-p_1},$$
thus, we find a super solution of (\ref{eq 1.1}) with $p=p_1$ and $\lambda\in(0,\lambda^*_{p_2})$,
and a minimal solution is derived for $\lambda\in(0,\lambda^*_{p_2})$. This means
\begin{equation}\label{increasing 1}
\lambda^*_{\gamma,p_2}\leq \lambda^*_{\gamma,p_1}\quad{\rm if}\ \ 0<p_1\leq p_2\leq p^*_\gamma.
\end{equation}

Thus for $p\in(0,p^*_\gamma]$ if $\gamma\in(0,1)$ or $p\in(0,1)$ if $\gamma=1$,  problem  (\ref{eq 1.1}) has the minimal solution $u_{p,\lambda}$ for  $\lambda\in(0,\lambda^*_p)$, 
$$ \lambda \mathbb{G}_\Omega [ a^{-p}] \leq  u_{p,\lambda} \leq a, $$
which implies (\ref{e 1.1-sub}) by Lemma \ref{lm 2.1}.

For $\lambda>\lambda^*_{\gamma,p}$, problem (\ref{eq 1.1}) has  no solution by the definition of $\lambda^*_{\gamma,p}$.\medskip

 {\it Finally, we show more properties of $\lambda^*_{\gamma,p}$. } $(iii)$ 
Now we choose  $p_1< p_2$ in $I_\gamma$, see (\ref{inter I}) for definition $I_\gamma$, and we recall  (\ref{increasing 1}), i.e.
 $\lambda^*_{\gamma,p_2}\leq \lambda^*_{\gamma,p_1}$ and  the map $\lambda\in(0,\lambda^*_{\gamma,p}] \mapsto u_{p,\lambda}$ is strictly increasing for $p\in I_\gamma$.  Then there exists $\delta_0\in(0,1)$ such that 
 $$0\leq a-u_{p_2,\lambda^*_{\gamma,p_2}}\leq a-u_{p_2,\lambda'}\leq \delta_0,$$
 where $\lambda'\in(0,\lambda^*_{\gamma,p_2})$ and $u_{p_2,\lambda'}$ is a positive classical solution of (\ref{eq 1.1})
 with $p=p_2$ and $\lambda=\lambda'$. 
 
From (\ref{4.6p3}), we have that $a-u_{p_1,\lambda^*_{p_1}}>0$ a.e. in $\Omega$ 
Let 
$$\theta_0=(a-u_{p_2,\lambda^*_{p_2}})^{p_1-p_2}\geq \delta_0^{p_1-p_2}>1$$
and then
\begin{eqnarray*} u_{p_2,\lambda^*_{\gamma,p_2}}=\lambda^*_{\gamma,p_2}\mathbb{G}_\Omega[(a-u_{p_2,\lambda^*_{\gamma,p_2}})^{-{p_2}}]\geq \theta_0\lambda^*_{\gamma,p_2}\mathbb{G}_\Omega[(a-u_{p_2,\lambda^*_{\gamma,p_2}})^{-{p_1}}].
\end{eqnarray*}
Therefore, $u_{p_2,\lambda^*_{\gamma,p_2}}$ is a super solution of (\ref{eq 1.1}) with $p=p_1$ and $\lambda=\theta_0\lambda^*_{\gamma,p_2}$,  from the proof of Theorem 1.1, we have that  (\ref{eq 1.1}) with $p=p_1$ and $\lambda=\theta_0\lambda^*_{\gamma,p_2}$ admits a minimal solution and 
$$\lambda^*_{\gamma,p_1}\geq \theta_0\lambda^*_{\gamma,p_2}>\lambda^*_{\gamma, p_2}.$$

\smallskip

$(iv)$ We recall that  $A_{\delta}=\{x\in\Omega:\, \rho(x)<\delta\}$ and $n\in\N$, $\eta_n$ is an $C^2$ function satisfying (\ref{101}), then 
\begin{eqnarray*}
|\int_\Omega u_{p,\lambda} (-\Delta \eta_n) dx|
 \leq  c_{11}n^2\int_{A_{\frac1n}} u_{p,\lambda} dx
 &\leq&
   c_{11}n^2\int_{A_{\frac1n}} a\, dx
   \\&\leq & c_{11}n^2\kappa \int_{A_{\frac1n}}\rho \,dx  \leq D_0,
\end{eqnarray*}
where   $|(-\Delta)\eta_n|\le c_7n^2$, $a(x)\leq \kappa \rho$,  $D_0$ is independent of $p$ and $u_{p,\lambda}$ is  the classical minimal solution of  (\ref{eq 1.1}) for any $\lambda<\lambda^*_{\gamma,p}$.  Then
$$\int_\Omega u_{p,\lambda} (-\Delta\eta_n) dx%=\int_\Omega \nabla u_{p,\lambda}\cdot \nabla{\eta_n} dx
=\int_\Omega (-\Delta u_{p,\lambda})\eta_n dx=\int_\Omega \frac{\lambda \eta_n}{(a-u_{p,\lambda})^p} dx.$$
Thus,
\begin{eqnarray*}
 D_0 &\geq &\int_\Omega   (-\Delta u_{p,\lambda})   dx\\
 &=& \lambda \int_\Omega \frac{1}{[a(x)-u_{p,\lambda}(x)]^p}dx
 \\& >&  \lambda \int_\Omega a^{-p }dx\geq \frac{\lambda}{\kappa^p}\int_{A_{\delta_0}} \rho^{-p}dx\geq c_{12}\frac{\lambda}{\kappa} (1-p)^{-1},
\end{eqnarray*}
where $c_{12}>0$ is independent of $p$. Thus
$$
\lambda<   \frac{\kappa  D_0}{c_{12} } (1-p)\to0\quad{\rm as}\ \ p\to1^-.
$$
So we have that $\lambda^*_{\gamma,p}\leq  \frac{\kappa  D_0}{c_{12} } (1-p)\to 0$ as $p\to1^-$.
\hfill $\Box$

  \medskip
  
\noindent{\bf Proof of Theorem \ref{teo 2}. }  By contradiction we assume that there exists some $\lambda>0$ such that
problem (\ref{eq 1.1}) has a solution $u_\lambda$ satisfying $0<u_\lambda<a$ in $\Omega$, then
\begin{equation}\label{2.5}
 \lambda \mathbb{G}_\Omega[ a^{-p}]\leq \lambda \mathbb{G}_\Omega[ (a-u_\lambda )^{-p}] = u \leq  a\quad{\rm in}\quad \Omega.
\end{equation}
On the other hand, by (\ref{1.1}), we have that
$$ u_\lambda\geq  \lambda \mathbb{G}_\Omega[ a^{-p}]\ge \frac{\lambda}{\kappa^p}\mathbb{G}_\Omega[\rho^{-p\gamma}]\geq   c_{13}\lambda \varrho_{_{2-p\gamma }}.$$
When $p>p^*_\gamma$ for $\gamma\in(0,1)$ or $p\geq 1$ for $\gamma=1$ or $p>0$ for $\gamma>1$  we have that 
$$\lim_{\rho(x)\to0^+}\mathbb{G}_\Omega[\rho^{-p\gamma}]\rho^{\gamma}=+\infty,$$
which, together with (\ref{2.5}),   deduces that
$$\lim_{x\in\Omega,  x\to\partial\Omega} a(x)\rho(x)^{-\gamma}=+\infty.$$
Then we obtain a contradiction from the assumption (\ref{1.1}).  
\hfill $\Box$

\smallskip

\subsection{ Extremal solutions }
We first prove the existence of the extremal solutions in weak sense by approximating the minimal solutions for $\lambda<\lambda^*_{\gamma,p}$.\medskip

 \noindent {\bf Proof of Theorem \ref{teo 5}. }  
 We recall that problem (\ref{eq 1.1}) with $\lambda\in(0,\lambda^*_{\gamma,p})$ admits a classical minimal solution  $ u_{p,\lambda}$.
For any $\beta\in(0,\gamma)$ and $n\in\N$, denote $\xi_n=\mathbb{G}_{\Omega}[\rho^{-1-\beta}]\eta_n$, where $\eta_n$ is defined by (\ref{101}),
 we observe that $\xi_n\in C_c^2(\Omega)$ and
\begin{eqnarray*}
\int_\Omega u_{p,\lambda} (-\Delta\xi_n) dx=\int_\Omega u_{p,\lambda} ( \rho^{-1-\beta}\eta_n-\nabla \mathbb{G}_{\Omega}[\rho^{-1-\beta}]\cdot\nabla\eta_n+(-\Delta)\eta_n\mathbb{G}_{\Omega}[\rho^{-1-\beta}])dx.
\end{eqnarray*}
Using Lemma \ref{lm 2.1} with $\tau=1-\beta\in(1-\gamma,1)$, we have that $\mathbb{G}_{\Omega}[\rho^{-1-\beta}]\le c_{14}\rho^{1-\beta}$ in $\Omega$.
Combining with the fact that $|\nabla \eta_n|\le c_5n$, $|(-\Delta)\eta_n|\le c_6n^2$
and  $0<u_{p,\lambda}(x)<a(x)\le \kappa \rho(x)^\gamma$ for $x\in\Omega$, we obtain that
\begin{eqnarray*}
% \nonumber to remove numbering (before each equation)
  \int_\Omega u_{p,\lambda} (-\Delta\xi_n) dx&\le& \kappa\int_\Omega  \rho^{\gamma-1-\beta} dx+c_{15}n\int_{A_{1/n}}  \rho^{\gamma-\beta}dx + c_{16}n^2\int_{A_{1/n}} \rho^{1-\beta+\gamma} dx \\
   &\le &   b_\beta,
\end{eqnarray*}
where $ c_{15}, c_{16},   b_\beta>0$ are independent on $n$ and $  b_\beta$ satisfies $  b_\beta\to+\infty$ as $\beta\to\gamma^-$. Thus,
$$\int_\Omega u_{p,\lambda} \rho^{-1-\beta} dx\leq b_\beta\quad{\rm and}\quad \int_\Omega \frac{ \lambda  \mathbb{G}_{\Omega}[\rho^{-1-\beta}]}{(a-u_{p,\lambda})^p} dx\le  b_\beta,$$
 then
\begin{equation}\label{4.6p}
\int_\Omega \frac{ \lambda  \rho^{1-\beta}}{(a-u_{p,\lambda})^p} dx\le  c_{17}  b_\beta .
\end{equation}
By Theorem \ref{teo 1}, we see that the mapping $\lambda\mapsto u_{p,\lambda}$ is increasing and
uniformly bounded by the function $a$, which is in $L^1(\Omega)$,  let
\begin{equation}\label{4.2}
 u_{p,\lambda^*_p}:=\lim_{\lambda\to\lambda^*_{\gamma,p}} u_{p,\lambda}\quad{\rm in}\ \ \bar\Omega.
\end{equation}
Then
$$u_{p,\lambda}\to u_{p,\lambda^*_{\gamma,p}}\ \ \  {\rm in}\ \ L^1(\Omega)\quad \ {\rm as}\ \lambda\to \lambda^*_{\gamma,p}$$
and then for any $\xi\in C_c^2(\Omega)$, we have that
\begin{equation}\label{4.6p1}
\int_\Omega u_{p,\lambda}(-\Delta\xi) dx\to \int_\Omega u_{p, \lambda^*_p}(-\Delta \xi) dx \ \ \ {\rm as} \ \ \lambda\to  \lambda^*_{\gamma,p}.
\end{equation}
Moreover,  we observe that the mapping $\lambda\mapsto  \frac{\lambda  }{(a-u_{p,\lambda})^p}$ is increasing and by (\ref{4.6p}),
$$\frac{\lambda  }{(a-u_{p,\lambda})^p}\to \frac{\lambda^*_{\gamma,p}  }{(a-u_{p,\lambda^*_{\gamma,p}})^p}\ \ \  {\rm a.e.\ in }\ \  \Omega\quad {\rm as}\ \lambda\to \lambda^*_{\gamma,p},$$
Then  (\ref{4.6p}) deduces that
$$\frac{\lambda  }{(a-u_{p,\lambda})^p}\to \frac{\lambda^*_{\gamma,p}  }{(a-u_{p,\lambda^*_{\gamma,p}})^p}\ \ \  {\rm in}\ \ L^1(\Omega,\,\rho^{1-\beta}dx)\quad {\rm as}\ \lambda\to \lambda^*_{\gamma,p}$$
and
\begin{equation}\label{4.6p3}
\int_\Omega \frac{ \lambda^*_{\gamma,p}  \rho^{1-\beta}}{(a-u_{p,\lambda^*_{\gamma,p}})^p} dx\le  c_{15} b_\beta ,
\end{equation}
thus, for any $\xi\in C_c^2(\Omega)$, we have that
\begin{equation}\label{4.6p2}
\int_\Omega \frac{\lambda \xi }{(a-u_{p,\lambda})^p} dx\to \int_\Omega \frac{\lambda^*_{\gamma,p} \xi }{(a-u_{p,\lambda^*_{\gamma,p}})^p} dx \ \ \ {\rm as} \ \ \lambda\to  \lambda^*_{\gamma,p}.
\end{equation}
Since $u_{p,\lambda}$ is a classical  solution of  (\ref{eq 1.1}),
$$
\int_\Omega u_{p,\lambda}(-\Delta\xi) dx=\int_\Omega \frac{\lambda \xi}{(a-u_{p,\lambda})^p} dx,\ \ \ \ \ \forall \, \xi\in C_c^2(\Omega),
$$
 passing to the limit as $\lambda\to\lambda^*_{\gamma,p}$, combining with (\ref{4.6p1}) and (\ref{4.6p2}),
then
 $u_{p,\lambda^*_{\gamma,p}}$ is a weak solution of (\ref{eq 1.1}) with $\lambda=\lambda^*_{\gamma,p}$.

By \cite[Theorem 2.6]{BV},   for any $\beta'\in(\beta,\gamma)$,  there exists $C_{\beta'}>0$ such that
$$
\norm{|\nabla u_{p,\lambda^*_{\gamma,p}}|}_{M^{\frac{N}{N- \beta'}}(\Omega) }\le C_{\beta'}\norm{(a-u_{p,\lambda^*_{\gamma,p}})^{-p}}_{L^1(\Omega, \, \rho^{1-\beta'}dx)},
$$
where ${M^{\frac{N}{N- \beta'}}(\Omega) }$ is the Marcinkiewicz space of
exponent $\frac{N}{N- \beta'}$.
Then, by (\ref{4.6p3}), we have that
$$\norm{|\nabla u_{p,\lambda^*_{\gamma,p}}|}_{L^{\frac{N}{N- \beta}}(\Omega)}\le   C_{1,\beta}(\lambda^*_{\gamma,p})^{-1}$$
and then
$$
 \norm{u_{p,\lambda^*_{\gamma,p}}}_{W^{1,\frac{N}{N- \beta}}(\Omega)}\le C_{2,\beta}
$$
for some $C_{i,\beta}>0$, $i=1,2$.
\hfill$\Box$  \medskip

 \setcounter{equation}{0}
\section{ Stability of the minimal solutions}

 \subsection{  A first result of    Stability }

In this subsection, we introduce the stability of the minimal solution $u_{p,\lambda}$ for problem (\ref{eq 1.1}) when the solution is comparable to $a$. Denote
\begin{equation}\label{5.0} 
\lambda_*=\sup\big\{\lambda>0:\ \liminf_{x\in\Omega,x\to\partial\Omega}(a(x)-u_{p,\lambda}(x))\rho(x)^{-\gamma}>0 \big\},
\end{equation}
 then  we have that $\lambda_\# \leq \lambda_* \leq\lambda^*_{\gamma,p}$, where $\lambda_\# $ is given in (\ref{def mu}).  Note that $\lambda_*$  depends on $p$.
 Then  for $\lambda\in (0,\lambda_*)$ there exists $C_\lambda\geq 1$ such that 
 $$\frac1{C_\lambda}\rho^{\gamma}\leq u_{p,\lambda}(x)\leq  C_\lambda \rho^{\gamma},\ \  \ \forall\,x\in \Omega.$$
 However it is difficult to show if $\lambda_*=\lambda^*_{\gamma,p}$.

By the definition of $\lambda_* $, we observe that for
any $\lambda\in (0,\lambda_*)$, there exists $\theta\in(0,1)$ such that $u_{p,\lambda}\le \theta  a$. Therefore,
by the fact that $a\ge \frac1\kappa \rho^\gamma$ in $\Omega$, we have that
\begin{equation}\label{5.1}
\frac{1}{(a-u_{p,\lambda})^{p+1}}\le \frac{1}{(1-\theta)^{p+1}a^{p+1}}\le (1-\theta)^{-(p+1)}\kappa^{p+1}\rho^{-(p+1)\gamma}\leq c_{16} \rho^{-2} \quad{\rm in}\ \, \Omega
\end{equation}
for $p\leq p^*_\gamma$ some $c_{16}>0$. 

Generally, we  consider the first eigenvalue of $-\Delta-\frac{p\lambda}{(a-u_{p,\lambda})^{p+1}}$ in $\cH_0^1(\Omega)$, that is,
\begin{equation}\label{eq eigen}
\mu_1(u,\lambda)=\inf_{\varphi\in \cH_0^1(\Omega)}\frac{\int_\Omega (|\nabla \varphi|^2-\frac{p\lambda \varphi^2}{(a-u)^{p+1}})dx}{\int_\Omega \varphi^2\,dx},
\end{equation}
which is well defined for $\lambda\in(0,\lambda^*_{\gamma,p})$, since $C^2_c(\Omega)\subset \cH_0^1(\Omega)$. For simplicity, we use the notation 
$$\mu_1(\lambda) =\mu_1(u_{p,\lambda},\lambda).$$
 It is well-known that  $u$ is stable if $\mu_1(u,\lambda)>0$  and semi-stable if $\mu_1(u,\lambda)\ge0$.

\begin{lemma}\label{lm 4.1}
Assume that $a\in   C^2_0(\Omega)$ satisfies (\ref{1.1}) with 
$\kappa\geq 1$ and $\gamma\in(0,1]$, $p\in \I_\gamma$.   If    $u$ is a  solution of (\ref{eq 1.1}) with $\lambda\in(0,\lambda^*_{\gamma,p}]$, $v $ is a super solution of (\ref{eq 1.1}) and $\mu_1(u,\lambda)$ is well-defined,
  then
$$u \le v \quad{\rm in}\ \, \Omega\quad {\rm if}\ \, \mu_1(u,\lambda)>0$$
and
$$u = v \quad{\rm in}\ \, \Omega\quad {\rm if}\ \, \mu_1(u,\lambda)=0.$$
\end{lemma}
{\bf Proof.}   For a given $\sigma\in[0,1]$, let $w_\sigma=u +\sigma(v -u )$
and then 
\begin{eqnarray*}
 \mathcal{F}(\sigma,x):&=&-\Delta w_\sigma-\frac{\lambda}{(a-w_\sigma)^p} \\[1mm]&\geq & \frac{\lambda(1-\sigma)}{(a-u)^p} +\frac{\lambda \sigma}{(a-v )^p}- \frac{\lambda}{(a-w_\sigma)^p}   \\[1mm]
   &\geq &  0, 
\end{eqnarray*}
since the function $t\mapsto \frac{\lambda}{(a(x)-t)^p}$ is convex on $(0,\, a(x))$. Note that
$ \mathcal{F}(0,x)=0$ and so $\partial_\sigma \mathcal{F}(0,x)\geq0$, i.e.
$$-\Delta(v -u )-\frac{p\lambda }{(a-u)^{p+1}}(v -u )\geq0\ \ {\rm in}\ \, \Omega,$$
 which, together with $v -u =0$ on $\partial\Omega$. Then the maximum principle implies that 
 if $\mu_1(u,\lambda)>0$, we have $v \geq u $ in $\Omega$. 
 
 If $\mu_1(u,\lambda)=0$, we have that
 $$-\Delta(v -u )-\frac{p\lambda }{(a-u )^{p+1}}(v -u )=0\ \ {\rm in}\ \, \Omega,$$ 
 which implies that $\partial_\sigma\mathcal{F}(0,x)=0$.  Since $\mathcal{F}(0,x)= 0$ and $\mathcal{F}(\sigma,x)\geq 0$ for $\sigma\in[0,1]$,   then  $\partial_{\sigma\sigma}\mathcal{F}(0,x)\geq0$ and  $$-\frac{p(p+1)\lambda }{(a-u )^{p+2}}(v -u)^2\geq0.$$
That holds only if
  $v =u $ in $\Omega$.
\hfill$\Box$\medskip

In order to show the well-definition of  $\mu_1(u,\lambda)$, we need the estimate of the decaying of $(a-u_{p,\lambda})^{-p-1}$, which is the key point to get the stability of minimal solutions. 

\begin{lemma}\label{lm kkk} 

Assume that $a\in C_0^2(\Omega) $ satisfies (\ref{1.1}) with 
$\kappa\geq 1$ and $\gamma\in(0,1]$, $p\in I_\gamma$. Let    $u_{p,\lambda}$ be the minimal  solution of (\ref{eq 1.1}) with $\lambda\in(0,\lambda_* )$.

Then   there exists $c_{16}>0$ such that  for any $\lambda\in(0,\lambda_* )$ 
\begin{equation}\label{5.1}
\frac{1}{(a-u_{p,\lambda})^{p+1}}  \leq c_{16}  (\lambda_*-\lambda )^{-p-1} \rho^{-2}  \quad{\rm in}\quad \Omega.
\end{equation}
\end{lemma}
{\bf Proof. } Since $\lambda\mapsto u_{p,\lambda}$ is increasing and verifies (\ref{e 1.1}) for fixed $p\in I_\gamma$,
then letting $\lambda'=\frac{\lambda_* +\lambda}{2}$  there holds that for $c_{17}>0$  
\begin{eqnarray*}
a-u_{p,\lambda}&\ge& u_{p,\lambda'}-u_{p,\lambda}  \\
 &\geq & (\lambda'-\lambda)\mathbb{G}_\Omega[a^{-p}]
  \geq c_{17}  (\lambda_* -\lambda ) \varrho_{_{2-p\gamma}}(x),\quad\forall\, x\in\Omega.
\end{eqnarray*}
 
For $ p\in(0,\, p^*_\gamma]$ with $\gamma\in(0,1)$,  Lemma \ref{lm kkk} implies that if $2-p\gamma\not=1$
\begin{eqnarray*}
\frac{1}{(a-u_{p,\lambda})^{p+1}}&\le& c_{17}^{-p-1}(\lambda_* -\lambda )^{-p-1} \varrho_{2-p\gamma}^{-p-1}
\\[1mm]&\leq& c_{17}^{-p-1}(\lambda_* -\lambda )^{-p-1} \rho^{-(p+1)\min\{1,2-p\gamma\}}  \quad{\rm in}\ \, \Omega,
\end{eqnarray*}
where we want to find $p$ such that 
$$-(p+1)\min\{1,2-p\gamma\} =-\min\{p+1,  (p+1)(2-p\gamma)\}\geq -2.$$
Note that $(p+1)(2-p\gamma)>2$ for $p\in(0,p^*_\gamma)$ and $(p+1)(2-p\gamma)=2$ for $p=p^*_\gamma$.  

If $p=\frac1\gamma\geq 1$, we have that 
\begin{eqnarray*}
\frac{1}{(a-u_{p,\lambda})^{p+1}}&\le& c_{17}^{-p-1} (\lambda_*-\lambda )^{-p-1} \varrho_{2-p\gamma}^{-p-1}
\\[1mm]&\leq& c_{17}^{-p-1}(\lambda_*-\lambda )^{-p-1} \rho^{-(p+1)} \ln\frac1{\rho}  \quad{\rm in}\ \, \Omega,
\end{eqnarray*}
where $ p+1 > 2$ for $\gamma\in(0,1)$.
Therefore, (\ref{5.1}) holds for  $p\in \I_\gamma.$ \hfill$\Box$\medskip

\begin{proposition}\label{pr 4.2}
Assume that $a\in   C^2_0(\Omega)$ satisfies (\ref{1.1}) with 
$\kappa\geq 1$ and $\gamma\in(0,1]$. Let    $u_{p,\lambda}$ be the minimal  solution of (\ref{eq 1.1}) with  $p\in \I_\gamma$ and $\lambda\in(0,\lambda_*)$,
then the minimal  solution $u_{p,\lambda}$ is stable.

\end{proposition}
{\bf Proof.}
  Notice that the mapping $\lambda\mapsto u_{p,\lambda}$ is strictly increasing and so is $\lambda\mapsto \frac{p\lambda}{(a-u_{p,\lambda})^{p+1}}$,
which implies that the mapping $\lambda\mapsto \mu_1(\lambda)$ is strictly decreasing
in $(0, \lambda_*)$. \smallskip

{\it Step 1:   $\mu_1(\lambda)>0$ and $u_{p,\lambda}$ is  stable   if $\lambda>0$ small.. }
It follows by \cite[Theorem 1]{MS} that there exists constant $c_{19}>0$ such that
\begin{equation}\label{4.1.1}
\int_\Omega \varphi^2\rho^{-2}dx\le c_{19}\int_{\Omega} |\nabla \varphi|^2 dx,\quad \forall \varphi\in \cH^1_0(\Omega).
\end{equation}
For $\lambda<\lambda^*_{\gamma,p}$, it follows from
 (\ref{5.1}) and (\ref{4.1.1}) that
$$\int_{\Omega}\frac{p\lambda \varphi^2}{(a-u_{p,\lambda})^{p+1}}dx\le p\lambda c_{16}\int_\Omega \varphi^2\rho^{-2}dx
\le p\lambda c_{16}c_{19} \int_{\Omega} |\nabla \varphi|^2 dx, \quad \forall \varphi\in \cH^1_0(\Omega).$$
For   $\lambda>0$ small such that $2\lambda c_{16}c_{19}< 1$, we obtain that
$$\int_{\Omega}\frac{p\lambda \varphi^2}{(a-u_{p,\lambda})^{p+1}}dx
<  \int_{\Omega} |\nabla \varphi|^2 dx, \quad \forall \varphi\in \cH_0^1(\Omega)\setminus\{0\}.$$

\smallskip

 {\it Step 2:   $u_\lambda$ is stable and $\mu_1(\lambda)>0$ for   $\lambda\in(0, \lambda_*)$. }    We have obtained that $\mu_1(\lambda)>0$ for $\lambda>0$ small and
 the mapping $\lambda\mapsto \mu_1(\lambda)$ is  strictly decreasing
in $(0, \lambda_*)$, so if there exists $\lambda_0\in(0,\  \lambda_*)$ such that $\mu_1(\lambda_0)=0$, then $\mu_1(\lambda)>0$ for $\lambda\in(0,\lambda_0)$.
Letting $\lambda_1\in (\lambda_0,\lambda_*)$, the minimal solution $u_{p,\lambda_1}$  is
a super solution of
$$
\left\{\arraycolsep=1pt
\begin{array}{lll}
 -\Delta    u = \frac{\lambda_0 }{(a-u)^p}\quad  &{\rm in}\quad\ \Omega,
 \\[2mm]
 \phantom{-\Delta   }
 u=0\quad &{\rm on}\quad   \partial \Omega
\end{array}
\right.
$$
and it infers from Lemma \ref{lm 4.1} that
$$u_{p,\lambda_1}=u_{p,\lambda_0},$$
 which contradicts that the mapping $\lambda\mapsto u_{p,\lambda}$ is strictly increasing  for   $\lambda\in(0, \lambda^*_{\gamma,p})$.
 Therefore, $\mu_1(\lambda)>0$  and  $u_{p,\lambda}$ is stable for $\lambda\in(0,  \lambda_*)$.
 \hfill$\Box$\medskip

\begin{lemma}\label{lm 5.1-=}
Assume that
$a\in C^2_0(\Omega)$ satisfies  (\ref{1.1})   with $\gamma\in(0,1), p=p^*_\gamma$. Then $\lambda^*_{\gamma,p}=\lambda_*$.

\end{lemma}
{\bf Proof.} For $\lambda\in(0,\lambda^*_{\gamma,p})$,
by (\ref{e 1.1}) and Lemma \ref{lm 2.1}, we have  that
 \begin{eqnarray*}
  a-u_{p,\lambda}& \ge& (\lambda^*_{\gamma,p}-\lambda)\mathbb{G}_\Omega[a^{-p}]\\[1mm]&\ge& c_{20}(\lambda^*_{\gamma,p}-\lambda)\rho^{2-p\gamma}
    \ge  c_{21}(\lambda^*_{\gamma,p}-\lambda)a,
 \end{eqnarray*}
 then there exists $\theta\in(0,1)$ such that $u_{p,\lambda}\le \theta a$ in $\Omega$. By   the definition of $\lambda_*$ that $\lambda_*=\lambda^*_{\gamma,p}$.
\hfill$\Box$

 \subsection{Full analysis of stability}

In order to who the stability of the minimal solution $u_{p,\lambda}$ for any $\lambda\leq \lambda^*_{\gamma,p}$, we consider the following equation. 
 \begin{equation}\label{eq 1.1-ep}
\left\{\arraycolsep=1pt
\begin{array}{lll}
\displaystyle -\Delta    u =  \frac{\lambda}{(a+\epsilon-u)^{p}}\quad  &{\rm in}\quad\ \Omega,
 \\[2.5mm]
 \phantom{- }
\displaystyle  0<u<a+\epsilon\quad &{\rm in}\quad\ \Omega,
 \\[1mm]
 \phantom{-\Delta   }
\displaystyle   u=0\quad &{\rm on}\quad   \partial \Omega,
\end{array}
\right.
\end{equation}
 where $\epsilon>0$. 
 \begin{lemma}\label{lm 3.1-ep}
 Assume that $\epsilon\in(0,1]$, $a\in C^2_0(\Omega)$ satisfies  (\ref{1.1})   with $\gamma\in(0,1],\, p\in\I_\gamma$. Then there exists $\lambda_\epsilon^*\geq \lambda^*_{\gamma,p}$ such that 
 (\ref{eq 1.1-ep}) has a minimal solution $w_{\epsilon,\lambda}$  for $\lambda\in(0,\lambda_\epsilon^*)$.
 
Furthermore, $(i)$ the mapping $\epsilon\in (0, 1]\to w_{\epsilon,\lambda}$ is  decreasing  and   for $\lambda\in (0,\lambda^*_{\gamma,p})$
$$w_{\epsilon,\lambda}<u_{p,\lambda}\quad{\rm in}\ \, \Omega. $$

$(ii)$   $w_{\epsilon,\lambda}$ is  stable for $\lambda<\lambda^*_{\epsilon}$ and semi-stable for $\lambda=\lambda^*_{\epsilon}$.

$(iii)$ There holds
$$\lim_{\epsilon\to0^+}w_{\epsilon,\lambda}=u_{p,\lambda} \quad{\rm in }\ C( \Omega)\ \, {\rm and}\ \ C^2(\mathcal{O})$$ 
 for any regular domain $\mathcal{O}$ such that  $\overline \mathcal{O}\subset \Omega$.
 \end{lemma}
 {\bf Proof. }    Let $v_{\epsilon,\lambda,0}\equiv0$ in $\bar\Omega$ and for $n\in\N$,
$$v_{\epsilon,\lambda,n}:=\lambda \mathbb{G}_\Omega[(a+\epsilon-v_{\epsilon,\lambda,n-1})^{-p}]>0,$$
Obviously, for $\lambda\in (0,\lambda^*_{\gamma,p})$, $u_{p,\lambda}$ is a super solution of 
(\ref{eq 1.1-ep}),  then as the proof of Theorem \ref{teo 1}, we obtain the limit of $\{v_{\epsilon,n}\}$, denoting $w_{\epsilon,\lambda}$,  is the minimal solution such that 
$$w_{\epsilon,\lambda}<u_{p,\lambda}\quad{\rm in}\ \, \Omega. $$

 Let us denote 
   $$\lambda^*_{\epsilon}=\sup\big\{\lambda>0:\,  (\ref{eq 1.1-ep}) \ {\rm has \  a \  minimal \  solution \ with\ such} \ \lambda \big\}.$$
Obviously, $\lambda^*_{\epsilon}\geq \lambda^*_{\gamma,p}>0$. 
For $0<\epsilon_1<\epsilon_2<1$, then 
$$w_{\epsilon_1,\lambda}= \lambda \mathbb{G}_\Omega[(a+\epsilon_1-w_{\epsilon_1,\lambda})^{-p}]>\lambda \mathbb{G}_\Omega[(a+\epsilon_2-w_{\epsilon_1,\lambda})^{-p}], $$
that is,  $w_{\epsilon_1,\lambda}$ is a super solution of (\ref{eq 1.1-ep}) with $\epsilon=\epsilon_2$ and then 
$$w_{\epsilon_1,\lambda}\geq w_{\epsilon_2,\lambda}. $$
We conclude that 
the mapping $\epsilon\in (0, 1]\to w_{\epsilon,\lambda}$ is    decreasing. Moreover, for $\lambda\in(0,\lambda^*_{\gamma,p})$
$$ v_{\epsilon,\lambda,n}<w_{\epsilon,\lambda}<u_{p,\lambda}<a\quad{\rm in}\ \, \Omega.$$

Denote $\mu_{\epsilon,1}(\lambda)$ the first eigenvalue of $-\Delta-\frac{p\lambda}{(a+\epsilon-w_{\epsilon,\lambda})^{p+1}}$ in $\cH_0^1(\Omega)$, that is,
$$\mu_{\epsilon,1}(\lambda)=\inf_{\varphi\in \cH_0^1(\Omega)\setminus\{0\}}\frac{\int_\Omega (|\nabla \varphi|^2-\frac{p\lambda \varphi^2}{(a+\epsilon-w_{\epsilon,\lambda})^{p+1}})dx}{\int_\Omega \varphi^2\,dx},$$
which is well defined for $\lambda\in(0,\lambda^*_{\epsilon})$.

Let $0<\epsilon_1<\epsilon_2<1$ and $\phi_{\epsilon_1}$ be the achieved function of $\mu_{\epsilon_1,1}(\lambda)$ with $\norm{\phi_{\epsilon_1}}_{L^2(\Omega)}=1$,
then we have that
\begin{eqnarray*}
 0<\mu_{\epsilon_2,1}(\lambda) -\mu_{\epsilon_1,1}(\lambda)  &\le& \int_\Omega \left(|\nabla \phi_{\epsilon_1 }|^2-\frac{p\lambda  \phi_{\epsilon_1}^2}{(a+\epsilon_2-w_{\epsilon_2,\lambda })^{p+1}}\right)dx 
 \\[2mm]
   && - \int_\Omega \left(|\nabla \phi_{\epsilon_1}|^2-\frac{p\lambda  \phi_{\epsilon_1}^2}{(a+\epsilon_1-w_{\epsilon_1,\lambda })^{p+1}}\right)dx 
   \\[2mm]&= &  \lambda p\int_\Omega \left(\frac{1}{(a+\epsilon_1-w_{\epsilon_1,\lambda })^{p+1}}-\frac{1}{(a+\epsilon_2-w_{\epsilon_2,\lambda })^{p+1}}\right)\phi_{\epsilon_1}^2dx
   \\[2mm]&\le &\lambda  p(p+1)(\epsilon_2-\epsilon_1)\int_\Omega\frac{ \phi_{\epsilon_1 }^2}{(a+\epsilon_1-w_{\epsilon_1, \lambda})^{p+2}} dx.
\end{eqnarray*}
Thus, the mapping $\epsilon\in(0,1]\mapsto \mu_{\epsilon,1}(\lambda)$ is increasing and locally Lipschitz continuous.

 Observe that 
 \begin{eqnarray*} 
 w_{\epsilon,\lambda}=\lambda \mathbb{G}_\Omega\big[ (a+\epsilon-w_{\epsilon,\lambda} )^{-p}\big] \leq \lambda \mathbb{G}_\Omega\big[ (a-u_{p,\lambda} )^{-p}\big],
\end{eqnarray*}
then for $\lambda<\lambda_*$, it follows Proposition \ref{pr 4.2} that
$$\int_{\Omega}\frac{p\lambda \varphi^2}{(a+\epsilon-w_{\epsilon,\lambda} )^{p+1}}dx\le \int_{\Omega}\frac{p\lambda \varphi^2}{(a-u_{p,\lambda})^{p+1}}dx
<  \int_{\Omega} |\nabla \varphi|^2 dx, \quad \forall \varphi\in \cH^1_0(\Omega).$$
 So we conclude that  $\mu_{\epsilon,1}(\lambda)>0$ and $w_{\epsilon,\lambda}$ is  stable   if $\lambda\in(0,\lambda_*)$.\smallskip

 Now we show  that  $w_{\epsilon,\lambda}$ is stable and $\mu_{\epsilon,1}(\lambda)>0$ for   $\lambda\in(0, \lambda^*_{\epsilon})$.     We have obtained that $\mu_{\epsilon,1}(\lambda)>0$ for $\lambda\in(0,\lambda_*)$,
 the mapping $\lambda\mapsto \mu_{\epsilon,1}(\lambda)$ is locally Lipschitz continuous, strictly decreasing in $(0, \lambda^*_{\epsilon})$, so if there exists $\lambda_0\in[\lambda_*,\  \lambda^*_{\epsilon})$ such that $\mu_{\epsilon,1}(\lambda_0)=0$, then $\mu_{\epsilon,1}(\lambda)>0$ for $\lambda\in(0,\lambda_0)$.
Letting $\lambda_1\in (\lambda_0,\lambda_*)$, the minimal solution $w_{\epsilon,\lambda_1}  $  is
a super solution of
$$
\left\{\arraycolsep=1pt
\begin{array}{lll}
 -\Delta    u = \frac{\lambda_0 }{(a+\epsilon-u)^p}\quad  &{\rm in}\quad\ \Omega,
 \\[2mm]
 \phantom{-\Delta   }
 u=0\quad &{\rm on}\quad   \partial \Omega
\end{array}
\right.
$$
and it infers from Lemma \ref{lm 4.1} that
$$w_{\epsilon,\lambda_1}  = w_{\epsilon,\lambda_0}\quad {\rm in}\ \, \Omega,$$
 which contradicts that the mapping $\lambda\mapsto w_{\epsilon,\lambda}$ is strictly increasing.
 Therefore, $\mu_{\epsilon,1}(\lambda)>0$  and  $w_{\epsilon,\lambda_1}$ is stable for $\lambda\in(0,  \lambda^*_{\epsilon})$.

 Finally,  from the property that the mapping $\epsilon\in (0, 1]\to w_{\epsilon,\lambda}$ is   decreasing and bounded by $u_{p,\lambda}\in C^2_0(\Omega)$,  then 
 the limit of $\{w_{\epsilon,\lambda}\}_{\epsilon}$ exists as $\epsilon\to0^+$ and it is a classical solution of (\ref{eq 1.1}). Since $u_{p,\lambda}$ is the minimal solution, we obtain that 
$$\lim_{\epsilon\to0^+}w_{\epsilon,\lambda}=u_{p,\lambda} \quad{\rm in }\ C( \Omega)\ \, {\rm and}\ \ C^2(\mathcal{O})$$ 
 for any compact domain $\mathcal{O}\subset \Omega$.

Then  in any compact domain $\mathcal{O}\subset \Omega$
\begin{equation}\label{5.0-con} 
\frac{p\lambda}{(a+\epsilon-w_{\epsilon,\lambda} )^{p+1}}- \frac{p\lambda}{(a-u_{p,\lambda} )^{p+1}}    \to0    \quad{\rm uniformly\ in}\ \mathcal{O}.
\end{equation}
 This completes the proof.  \hfill$\Box$\medskip

 \begin{remark}
 $(I)$ For the existence of minimal solution of (\ref{eq 1.1-ep}), the restriction of $p$ could be enlarged to the interval $(0,+\infty)$.  
 
  $(II)$ The pull-in voltage $\lambda^*_\epsilon$ could be proved to be increasing with respect to $\epsilon$. 
 
 \end{remark}

\noindent{\bf Proof of Theorem \ref{teo 3}. }  From Proposition \ref{pr 4.2}, we have that 
$u_{p,\lambda}$  is  stable for $p\in I_\gamma$ and $\lambda\in(0,\lambda_*)$. 

  From Lemma \ref{lm 5.1-=},  $\lambda_* =\lambda^*_{\gamma,p}$ for $\gamma\in(0,1)$ and $p=p^*_\gamma$.  Now we only have to deal with the case  that  
$$\lambda_*<\lambda^*_{\gamma,p}$$
which may happen  when $\gamma\in(0,1]$ and  $p\in(0,p^*_\gamma).$

\smallskip

{\it Step 1: we  show that $u_{p,\lambda}$  is  stable
for $p\in I_\gamma$ and $\lambda\in[\lambda_*, \lambda^*_{\gamma,p})$.  }\smallskip

From Lemma \ref{lm 3.1-ep},   we have that  $w_{\epsilon, \lambda}$ is stable for $\lambda\in[\lambda_*, \lambda^*_{\gamma,p})$ and $\epsilon>0$, that is,
$$
\int_\Omega\frac{p\lambda\varphi^2}{(a+\epsilon-w_{\epsilon, \lambda})^{p+1}}dx<\int_\Omega|\nabla \varphi|^2dx,\quad \ \forall \varphi\in \cH^1_0(\Omega)\setminus\{0\}.
$$
Let $\varphi=\mathbb{G}_{\Omega}[1]$, we have that
$$\int_\Omega\frac{ \rho^2}{(a+\epsilon-w_{\epsilon, \lambda})^{p+1}}dx<c_{22}\lambda^{-1},$$
where $c_{22}>0$ is independent of $\epsilon$. Since $\epsilon\in(0,1]\to w_{\epsilon, \lambda}$ is  decreasing and $$\lim_{\epsilon\to0^+}w_{\epsilon,\lambda}=u_{p,\lambda} \quad{\rm in }\ C( \Omega)\ \, {\rm and}\ \ C^2(\mathcal{O})$$ 
 for any compact domain $\mathcal{O}\subset \Omega$.
 then
  for any $\varphi\in C_c^2(\Omega)$,
$$\lim_{\epsilon\to0^+ }\int_\Omega\frac{\lambda\varphi^2}{(a+\epsilon-w_{\epsilon, \lambda})^{p+1}}dx=\int_\Omega\frac{\lambda  \varphi^2}{(a-u_{p,\lambda})^{p+1}}dx$$
Therefore,
$$\int_\Omega\frac{p\lambda \varphi^2}{(a-u_{p,\lambda })^{p+1}}dx\le \int_\Omega|\nabla \varphi|^2dx,\quad \forall \varphi\in C_c^2(\Omega)$$
by the fact that $C_c^2(\Omega)$ is dense in $\cH^1_0(\Omega)$,
then $u_{p,\lambda}$ is semi-stable.  So we obtain that  $\mu_1(\lambda)\geq0$. 

Now we prove $\mu_1(\lambda)>0$ for $[\lambda_*,\lambda^*_{\gamma,p})$.  If not, there
exists $\lambda_0\in(\lambda_*,\, \lambda^*_{\gamma,p})$ such that  $\mu_1(\lambda_0)=0$.

Note that  $\lambda\in [\lambda_*,\lambda^*_{\gamma,p})\to u_{p,\lambda}$
 is strictly increasing, take $\bar \lambda\in(\lambda_0, \lambda^*_{\gamma,p})$
 $u_{p,\lambda_0}<u_{p,\bar \lambda}$ and $u_{p,\bar \lambda}$ is a super solution of (\ref{eq 1.1})
 with $\lambda=  \lambda_0$. However, by the assumption that  $\mu_1(\lambda_0)=0$,  it follows by Lemma \ref{lm 4.1} that  $u_{p,\lambda_0}=u_{p,\bar \lambda }$, which is impossible.
 Thus, we conclude that $\mu_1(\lambda)>0$ and $u_{p,\lambda}$ is  stable for $\lambda\in [\lambda_*,\lambda^*_{\gamma,p})$.

{\it Step 2: we  show that $u_{p,\lambda}$  is semi-stable
for $p\in I_\gamma$ and $\lambda=\lambda^*_{\gamma,p}$.  }\smallskip

From {\it Step 1},  we have that  $u_{p,\lambda}$ is stable for $\lambda\in(0, \lambda^*_{\gamma,p})$, that is,
$$
\int_\Omega\frac{p\lambda\varphi^2}{(a-u_{p,\lambda})^{p+1}}dx<\int_\Omega|\nabla \varphi|^2dx,\quad \ \forall \varphi\in \cH^1_0(\Omega)\setminus\{0\}.
$$
Let $\varphi=\mathbb{G}_{\Omega}[1]$, we have that
$$\int_\Omega\frac{ \rho^2}{(a-u_{p,\lambda})^{p+1}}dx<c_{39}\lambda^{-1}.$$
Since $\lambda\in(0,\lambda^*_{\gamma,p})\to u_{p,\lambda}$
 is strictly increasing and continuous,  so is the mapping $\lambda\mapsto \frac{ \rho^2}{(a-u_{p,\lambda})^{p+1}}$, then $\frac{ \rho^2}{(a-u_{p,\lambda})^{p+1}}$ is bounded in $L^1(\Omega)$,
then
$$\frac{ \rho^2}{(a-u_{p,\lambda})^{p+1}}\to \frac{ \rho^2}{(a-u_{p,\lambda^*_{\gamma,p} })^{p+1}}\quad{\rm as}\quad \lambda\to \lambda^*_{\gamma,p} \quad{\rm in}\quad L^1(\Omega)$$
and for any $\varphi\in C_c^2(\Omega)$,
$$\lim_{\lambda\to \lambda^*_{\gamma,p} }\int_\Omega\frac{\lambda\varphi^2}{(a-u_{p,\lambda})^{p+1}}dx=\int_\Omega\frac{\lambda^*_{\gamma,p} \varphi^2}{(a-u_{p,\lambda^*_{\gamma,p}})^{p+1}}dx.$$
Therefore,
$$\int_\Omega\frac{p\lambda^*_{\gamma,p}\varphi^2}{(a-u_{p,\lambda^*_{\gamma,p}})^{p+1}}dx\le \int_\Omega|\nabla \varphi|^2dx,\quad \forall \varphi\in C_c^2(\Omega),$$
by the fact that $C_c^2(\Omega)$ is dense in $\cH^1_0(\Omega)$,
then $u_{p,\lambda^*_{\gamma,p}}$ is semi-stable and  $\mu_1(\lambda^*_{\gamma,p})\geq0$. 
\hfill$\Box$

\subsection{Continuity of the map $\lambda\mapsto u_{p,\lambda}$}
 
\begin{proposition}\label{pr 4.2-con}
Assume that $a\in   C^2_0(\Omega)$ satisfies (\ref{1.1}) with 
$\kappa\geq 1$ and $\gamma\in(0,1]$.  Let  $u_{p,\lambda}$ be the minimal  solution of (\ref{eq 1.1}) with  $p\in \I_\gamma$ and $\lambda\in(0,\lambda^*_{\gamma,p})$,
then   
$\lambda\in (0,\lambda^*_{\gamma,p})\mapsto u_{p,\lambda}$ is  continuous in $L^\infty(\Omega)$.
\end{proposition}
{\bf Proof. } {\it Now we show that 
the map $\lambda\in(0,\lambda^*_{\gamma,p})\mapsto u_{p,\lambda}$ is left continuous. }
Indeed, if not, let $u_{p,\lambda}$ is not left continuous at $\lambda_0\in(0,\lambda^*_{\gamma,p})$. From the increasing monotonicity $\lambda\in(0,\lambda^*_{\gamma,p})\mapsto u_{p,\lambda}$ and the upper bound $u_{p,\lambda_0}$, the limit of minimal solutions $\{u_{p,\lambda}\}_{\lambda\in(0,\lambda_0)}$ exists as $\lambda\to\lambda_0^-$,
$$\bar u_{p,\lambda_0}:=\lim_{\lambda\to\lambda_0^-}u_{p,\lambda}\quad{\rm in}\ \, \Omega,  $$
which is a classical solution of (\ref{eq 1.1}) with $\lambda=\lambda_0$ by the regularity and 
$$\bar u_{p,\lambda_0}\leq u_{p,\lambda_0}\quad{\rm in}\ \, \Omega.$$
Since $u_{p,\lambda_0}$ is a minimal solution, then both have to coincide. This is a contradiction and the argument is proved.

On the other hand, the lower bound $u_{p,\lambda_0}$ and upper bounded $u_{p,\lambda_0+\epsilon_0}$ with $\epsilon_0\leq \frac{\lambda^*_{\gamma,p}-\lambda_0}{2}$, the limit of minimal solutions $\{u_{p,\lambda}\}_{\lambda\in\big(\lambda_0,\frac{\lambda_0+\lambda^*_{\gamma,p}}{2}\big)}$ exists as $\lambda\to\lambda_0^+$,
$$\tilde u_{p,\lambda_0}:=\lim_{\lambda\to\lambda_0^+}u_{p,\lambda}\quad{\rm in}\ \, \Omega,  $$
which is a classical solution of (\ref{eq 1.1}) with $\lambda=\lambda_0$ by the regularity and 
$$\tilde u_{p,\lambda_0}\geq u_{p,\lambda_0}\quad{\rm in}\ \, \Omega.$$
Since $u_{p,\lambda}$ is a minimal solution and stable for $\lambda\in(\lambda_0,\frac{\lambda_0+\lambda^*_{\gamma,p}}{2})$, 
denoting 
$ \mu_1(\tilde u_{p,\lambda_0},\lambda)$  the eigenvalue of (\ref{eq eigen}) with $u=\tilde  u_{p,\lambda_0}$, then we have that $\mu_1(\tilde u_{p,\lambda_0},\lambda)> 0$, which implies by the strictly decreasing $\lambda\to\mu_1(\lambda)$. 
 Then we apply Lemma \ref{lm 4.1}
to obtain that $$\tilde u_{p,\lambda_0}\leq u_{p,\lambda_0}\quad{\rm in}\ \, \Omega.$$

Therefore, 
$$\lim_{\lambda\to\lambda_0}\|u_{p,\lambda_0}-u_{p,\lambda}\|_{L^\infty(\Omega)}=0$$
and we have that the map $\lambda\in(0,\lambda^*_{\gamma,p})\mapsto u_{p,\lambda}$ is continuous. %Moreover, 
%\begin{eqnarray*}
%u_{p,\lambda}-u_{p,\lambda_0}&=&\lambda  \mathbb{G}_\Omega[ (a-u_{p,\lambda})^{-p}]-\lambda_0  \mathbb{G}_\Omega[ (a-u_{p,\lambda_0})^{-p}]
%\\[2mm]&\leq & |\lambda-\lambda_0|\mathbb{G}_\Omega[ (a-u_{p,\lambda})^{-p}]+
%\lambda_0\mathbb{G}_\Omega\big[\big|(a-u_{p,\lambda})^{-p}-(a-u_{p,\lambda_0})^{-p}\big|\big]
%\\[2mm]&\leq &\frac1\lambda  |\lambda-\lambda_0| u_{p,\lambda}+p\lambda_0\mathbb{G}_\Omega\big[(a-u_{p,\max\{\lambda,\lambda_0\}})^{-p-1} |u_{p,\lambda_0}-u_{p,\lambda}|  \big]
%\end{eqnarray*}
\hfill$\Box$

\begin{proposition}\label{pr 4.2-con lip}
Assume that $a\in   C^2_0(\Omega)$ satisfies (\ref{1.1}) with 
$\kappa\geq 1$ and $\gamma\in(0,1]$.  Let  $u_{p,\lambda}$ be the minimal  solution of (\ref{eq 1.1}) with  $p\in \I_\gamma$ and $\lambda\in(0,\lambda_*)$,
then   there exists $\bar \lambda\in(0,\lambda_*]$ such that 
$\lambda\in (0,\bar \lambda )\mapsto u_{p,\lambda}$ is Lipschitz continuous.
\end{proposition}
{\bf Proof. } Now we show that  $\lambda\mapsto u_{p,\lambda}$ is Lipschitz continuous when $\lambda>0$ small enough. As the proof of Theorem \ref{teo 1}, the minimal solution $u_{p,\lambda}$ is approximated by the sequence $v_{\lambda,n}$ defined by 
$$v_{\lambda,1}=\lambda \mathbb{G}_\Omega[ a^{-p}]>0,\quad v_{\lambda,n}=\lambda \mathbb{G}_\Omega[ (a-v_{\lambda,n-1})^{-p}]\quad {\rm for}\ \,  n\in\N.$$
Now let $0<\lambda_1<\lambda_2\leq \bar \lambda\leq \lambda_*$, where for some   $$\bar\lambda\leq \frac1{2 p}\min\big\{1,2^{-p-1}c_0^{-1} \big\}$$  small such that 
$$u_{p,\bar \lambda}\leq \frac12 a$$
then  
$$v_{\lambda_1,n}<v_{\lambda_2,n}<u_{p,\lambda_2},\qquad    v_{\lambda_2,1}-v_{\lambda_1,1}=(\lambda_2-\lambda_1)\mathbb{G}_\Omega[ a^{-p}]\leq {\bf c}_0 (\lambda_2-\lambda_1)a  $$
and
\begin{eqnarray*}
0<v_{\lambda_2,2}-v_{\lambda_1,2}&=&(\lambda_2-\lambda_1) \mathbb{G}_\Omega[ (a-v_{\lambda_2,1})^{-p}]+\lambda_1\mathbb{G}_\Omega[ (a-v_{\lambda_2,1})^{-p}- (a-v_{\lambda_1,1})^{-p}]
\\[2mm]&\leq& 2^p(\lambda_2-\lambda_1) \mathbb{G}_\Omega[ a^{-p}]+p \lambda_1 \mathbb{G}_\Omega[ (a-v_{\lambda_2,1})^{-p-1}(v_{\lambda_2,1}-v_{\lambda_1,1})]
\\[2mm]&\leq& 2^p{\bf c}_0 (\lambda_2-\lambda_1) a+2^{p+1}p \lambda_1 (\lambda_2-\lambda_1)\mathbb{G}_\Omega[a ^{-p}]
\\[2mm]&\leq&  {\bf c}_1 (\lambda_2-\lambda_1)a,
 \end{eqnarray*}
where ${\bf c}_1={\bf c}_0(2^p+2^{p+1}p\lambda_1)\leq 2^{p+1}c_0$
since
 $\bar\lambda\leq  2^{-1}p^{-1}. $ 
 
Inductively, we assume that 
$$v_{\lambda_2,n+1}-v_{\lambda_1,n+1} \leq {\bf c}_n (\lambda_2-\lambda_1)a,$$
where we set  for $n\geq 2$
$${\bf c}_n\leq 2^{p+1}{\bf c}_0.$$

Then we have that 
\begin{eqnarray*}
0<v_{\lambda_2,n+1}-v_{\lambda_1,n+1}&=&(\lambda_2-\lambda_1) \mathbb{G}_\Omega[ (a-v_{\lambda_2,n})^{-p}]+\lambda_1\mathbb{G}_\Omega[ (a-v_{\lambda_2,n})^{-p}- (a-v_{\lambda_1,n})^{-p}]
\\[2mm]&\leq& 2^p(\lambda_2-\lambda_1) \mathbb{G}_\Omega[ a^{-p}]+p \lambda_1 \mathbb{G}_\Omega[ (a-v_{\lambda_2,n})^{-p-1}(v_{\lambda_2,n}-v_{\lambda_1,n})]
\\[2mm]&\leq& 2^p{\bf c}_0 (\lambda_2-\lambda_1) a+2^{p+1}p {\bf c}_n \lambda_1 (\lambda_2-\lambda_1)\mathbb{G}_\Omega[a ^{-p}]
\\[2mm]&\leq&  {\bf c}_0(2^p+2^{p+1}p\lambda_1{\bf c}_n) (\lambda_2-\lambda_1)a
 \end{eqnarray*}
That is,
$$0<v_{\lambda_2,n+1}-v_{\lambda_1,n+1}\leq {\bf c}_n (\lambda_2-\lambda_1)a, $$
where 
\begin{eqnarray*} 
{\bf c}_{n+1}&=&{\bf c}_0(2^p+2^{p+1}p\lambda_1 {\bf c}_{n})
\\[2mm]&\leq& {\bf c}_0(2^p+2^{p+1}p\, \bar \lambda\, 2^{p+1}c_0 ) 
 \leq  2^{p+1}{\bf c}_0,
\end{eqnarray*}
since
$$\bar\lambda\leq  2^{-p-2}p^{-1}{\bf c}_0^{-1}.$$
Therefore,   for $0<\lambda_1<\lambda_2 \leq \bar \lambda$ there holds
\begin{equation}\label{5.0-con1}  
0<u_{p,\lambda_2}-u_{p,\lambda_1}\leq  2^{p+1} {\bf c}_0 (\lambda_2-\lambda_1)a\quad{\rm in}\ \, \Omega.
\end{equation}
The proof ends. \hfill$\Box$\medskip

 \begin{corollary}\label{cr 4.2-con lip}
Assume that $a\in   C^2_0(\Omega)$ satisfies (\ref{1.1}) with 
$\kappa\geq 1$ and $\gamma\in(0,1]$. Let $\bar \lambda$ be given in Proposition \ref{pr 4.2-con lip}. 
Then %$(i)$ the map  $\lambda\in(0,  \lambda^*_{\gamma,p})\mapsto \mu_1(\lambda)$ is continuous;\smallskip
%$(ii)$ 
 the map  $\lambda\in (0,\bar \lambda )\mapsto \mu_1(\lambda)$ is Lipschitz continuous.
\end{corollary}
{\bf Proof. } Let $0<\lambda_1<\lambda_2\leq  \bar\lambda$ and $\phi_{\lambda_2}$ be the achieved function of $\mu_1(\lambda_2)$ with $\norm{\phi_{\lambda_2}}_{L^2(\Omega)}=1$,
then we have that
\begin{eqnarray*}
0< \mu_1(\lambda_1)-\mu_1(\lambda_2) &\le& \int_\Omega \left(|\nabla \phi_{\lambda_2}|^2-\frac{p\lambda_1 \phi_{\lambda_2}^2}{(a-u_{p,\lambda_1})^{p+1}}\right)dx  \\
   && -\int_\Omega \left(|\nabla \phi_{\lambda_2}|^2-\frac{p\lambda_2 \phi_{\lambda_2}^2}{(a-u_{p,\lambda_2})^{p+1}}\right)dx
   \\&\le &p(\lambda_2-\lambda_1)\int_\Omega\frac{ \phi_{\lambda_2}^2}{(a-u_{p,\lambda_2})^{p+1}} dx
   \\&&+p\lambda_1\int_\Omega\Big( \frac{1}{(a-u_{p,\lambda_2})^{p+1}} -\frac{1}{(a-u_{p,\lambda_1})^{p+1}}\Big)\phi_{\lambda_2}^2dx
   \\&\le & 2^{p+1} (\lambda_2-\lambda_1)\int_\Omega\frac{ \phi_{\lambda_2}^2}{ a^{p+1}} dx
  +p\lambda_1  (p+1) \int_\Omega \frac{u_{p,\lambda_2}-u_{p,\lambda_1}}{(a-u_{p,\lambda_2})^{p+2}} \phi_{\lambda_2}^2dx.
\end{eqnarray*}
%From Proposition \ref{pr 4.2-con},    $u_{p,\lambda_2}-u_{p,\lambda_1}=o(1) a$ in $\Omega$, where $o(1)\to0$ as
%$\lambda_2-\lambda_1\to0$, we have that 
%$$\lim_{\lambda_2-\lambda_1\to0}\Big(\mu_1(\lambda_1)-\mu_1(\lambda_2)\Big)=0.$$
From (\ref{5.0-con1}), we obtain that  
\begin{eqnarray*}  
 (p+1) \int_\Omega \frac{u_{p,\lambda_2}-u_{p,\lambda_1}}{(a-u_{p,\lambda_2})^{p+2}} \phi_{\lambda_2}^2dx  \le   (p+1) 2^{2p+3} {\bf c}_0 (\lambda_2-\lambda_1) \int_\Omega\frac{ \phi_{\lambda_2}^2}{ a^{p+1}} dx,
\end{eqnarray*}
then the mapping $\lambda\in(0,\bar \lambda)\mapsto \mu_1(\lambda)$ is locally Lipschitz continuous.\hfill$\Box$ \medskip

\begin{remark}
Because of the boundary decay of $a$, it is difficulty to obtain the Lipschitz continuity of $\lambda\mapsto u_{p,\lambda}$ for $\lambda>\bar\lambda$. 

\end{remark}

\setcounter{equation}{0}
\section{Regularity of extremal solutions }

% \subsection{Low dimension}

%We next improve the regularity of $u_{\lambda^*}$ and prove when $1\le N\le 7$, the extremal
%solution $u_{\lambda^*}$ is a classical solution of (\ref{eq 1.1}) with $\lambda=\lambda^*$. 
In order to improve the regularity of extremal solutions, 
we need the following Lemma.

\begin{lemma}\label{lm 4.2}
Assume that $a\in   C^2_0(\Omega)$ satisfies (\ref{1.1}) with 
$\kappa\geq 1$, $\gamma\in(0,1]$, $p\in\I_\gamma$ and $\lambda\in(0,\lambda^*_{\gamma,p}]$.  Let $u$ be a weak solution of (\ref{eq 1.1}) satisfying   such that
 \begin{equation}\label{4.1.3}
 \Big\|\frac1{(a-u)^p}\Big\|_{L^{q_\sigma}(\Omega_r)}<+\infty,
 \end{equation}
 where $q_\sigma=\frac{N}{p(2-\sigma)}$ with $\sigma\in[1,2)$ and   
 $$ \Omega_r:=\{y\in\Omega:\,\rho(y)>r\}\quad{\rm with}\ \, r>0.$$
 Then  
\begin{equation}\label{4.2.1}
 \inf_{x\in \Omega_{3r}}\big(a(x)-u(x)\big)>0
\end{equation}
and 
$u$ is $C^2$ in $\Omega_{3r}$. 

\end{lemma}
{\bf Proof.}  By (\ref{4.1.3}), we have that
$$\frac1{(a-u)^p}\in L^{\frac{N}{2-\sigma}}(\Omega_r)$$
and then $u\in W^{2, \frac{N}{2-\sigma}}(\Omega_{2r})$ and by Sobolev's Theorem, we
deduce that $u\in C^{\sigma}(\Omega_{3r})$. Here $C^{\sigma}=C^{0,1}$ if $\sigma=1$, $C^{\sigma}=C^{1,\sigma-1}$ if $\sigma\in(1,2)$. 

We next show that $u<a$ in $\Omega_{3r}$. Indeed, if not, there exists $x_0\in \Omega_{3r}$ such that $u(x_0)=a(x_0)$, which means that $u$ contacts $a$ at $x_0$.   Thus,
\begin{eqnarray*}
 |a(x)-u(x)|&\le &|a(x)-a(x_0)|+|u(x)-u(x_0)|+|u(x_0)-a(x_0)|\\[1mm] &=&
|a(x)-a(x_0)|+|u(x)-u(x_0)|
\\[1mm]&\le& c_{22}|x-x_0|^{\sigma},
\end{eqnarray*}
then for $\sigma\in[1,2)$ one has that 
$$\int_{\Omega_{3r}}\frac1{(a-u)^{pq_\sigma}}dx\geq\int_{\Omega_{3r}} |x-x_0|^{-\frac{\sigma}{2-\sigma}N}dx  \ge c_{23} \int_{\Omega_{3r}} |x-x_0|^{-N}dx=\infty,$$
which  contradicts (\ref{4.1.3}). Therefore, (\ref{4.2.1}) holds.
 \hfill$\Box$ \medskip

%While $p^*_\gamma=\frac{2}{\gamma}-1$, so we have that $p^*_\gamma\geq 1>\frac13>p^\#_\gamma$. 

\begin{lemma}\label{lm 8.1}
Assume that 
$a\in  C^2_0(\Omega)$    satisfies (\ref{1.1}) with 
$\kappa\geq 1$ and
 $\gamma\in(0,\, 1]$.
 %and
 %$$\int_\Omega |\Delta a|dx <+\infty.$$
  Let $p\in \I_\gamma$, 
 and  $u_{p,\lambda^*_{\gamma,p}}$ be the extremal solution of (\ref{eq 1.1}).
 Then 
$u_{p,\lambda^*_{\gamma,p}}\in W^{2,q}_{loc}(\Omega)$  for any  $1\leq q<f_0(p)$.
\end{lemma}
{\bf Proof.}   For $\lambda\in(0,\lambda^*_{\gamma,p})$, we know that $u_{p,\lambda}$ is stable, then
\begin{equation}\label{4.3.1}
\int_{\Omega}\frac{p\lambda \varphi^2}{(a-u_{p,\lambda})^{p+1}}dx< \int_{\Omega} |\nabla \varphi|^2 dx, \quad \forall \varphi\in \cH^1_0({\Omega})\setminus\{0\}.
\end{equation}

Let us denote
$$
\varphi_\theta= \left\{\arraycolsep=1pt
\begin{array}{lll}
 (a -u_{p,\lambda})^\theta \eta_r^\tau \quad  &{\rm in}\ \  \Omega,
 \\[1.5mm]
 \phantom{}
0\quad  &{\rm on}\ \, \partial \Omega,
\end{array}
\right.
$$
where $\theta\in\big(-p-\sqrt{p+p^2},0\big)$, $\tau>1$ and 
   $\eta_r:\Omega \to[0,1]$ is a function in $C_c^1(\Omega)$ such that $\eta_r=1$ in $\Omega_{2r} $ and $\eta_r=0$ in $\Omega\setminus \Omega_r$. We observe that $\varphi_\theta\in \cH_0^1(\Omega)$.
It follows by (\ref{4.3.1}) with $\varphi=\varphi_\theta$  that
\begin{eqnarray}
  p\lambda  \int_{\Omega_r}\frac{ (a-u_{p,\lambda})^{2\theta} \eta_r^{2\tau}}{(a-u_{p,\lambda})^{p+1}}dx
 &\le &
 \int_{\Omega}|\nabla \big((a-u_{p,\lambda})^\theta\eta_r^\tau\big)|^2dx\nonumber \\
   &\leq &  (\theta^2+\delta)\int_{\Omega}(a-u_{p,\lambda})^{2\theta-2}  |\nabla (a-u_{p,\lambda})|^2 \eta_r^{2\tau} dx\nonumber
   \\&&+(1+\frac1\delta)\int_{\Omega}(a-u_{p,\lambda})^{2\theta} |\nabla (\eta_r^\tau)|^2dx,\label{4.3.2}
\end{eqnarray}
where $\delta\in(0,\frac{1}4\theta^2)$ small enough.
Since $u_{p,\lambda}$ is the minimal solution of (\ref{eq 1.1}), then
\begin{equation}\label{4.3.3}
  -\Delta (a-u_{p,\lambda})=-\Delta a-\frac{\lambda}{(a-u_{p,\lambda})^p}\quad {\rm in} \ \  \Omega.
\end{equation}
Multiplying  (\ref{4.3.3})  by $\frac{\theta^2}{1-2\theta}  (a-u_{p,\lambda})^{2\theta-1}\eta_r^{2\tau}$ and then integration by parts yields that
\begin{eqnarray*}
&&\frac{\theta^2}{1-2\theta}\int_{\Omega}\Big( \Delta a+\frac{\lambda}{(a-u_{p,\lambda})^p}\Big) (a-u_{p,\lambda})^{2\theta-1} \eta_r^{2\tau}dx
\\[2mm]  & =&-\frac{\theta^2}{1-2\theta}\int_{\Omega} \Delta (a-u_{p,\lambda}) (a-u_{p,\lambda})^{2\theta-1} \eta_r^{2\tau}dx
  \\[2mm]  & =&-\frac{\theta^2}{1-2\theta}\int_{\Omega}\nabla (a-u_{p,\lambda})\cdot\nabla\left((a-u_{p,\lambda})^{2\theta-1}\right)\eta_r^{2\tau} dx
  \\[2mm]&& - \frac{\theta^2}{1-2\theta}\int_{\Omega}  (a-u_{p,\lambda})^{2\theta-1} \nabla (a-u_{p,\lambda})\cdot \nabla (\eta_r^{2\tau}) dx
 \\[2mm] & \geq &  \theta^2 \int_{\Omega}(a-u_{p,\lambda})^{2\theta-2}|\nabla(a-u_{p,\lambda})|^2\eta_r^{2\tau} dx
  - \frac{\theta^2}{1-2\theta}\int_{\Omega}  (a-u_{p,\lambda})^{2\theta-1}  |\nabla (a-u_{p,\lambda})| |\nabla (\eta_r^{2\tau})| dx
   \\[2mm]  & >& \theta^2 \int_{\Omega}(a-u_{p,\lambda})^{2\theta-2}|\nabla(a-u_{p,\lambda})|^2\eta_r^{2\tau} dx 
    \\[2mm]&& - \frac{2 \theta^2}{1-2\theta}  \Big(\int_{\Omega}  (a-u_{p,\lambda})^{2\theta-2}  |\nabla (a-u_{p,\lambda})|^2 \eta_r^{2\tau} dx\Big)^{\frac12}
     \Big(\int_{\Omega}  (a-u_{p,\lambda})^{2\theta}|\nabla \eta_r|^2dx\Big)^{\frac12}
    \\[2mm]  & >&(\theta^2-\delta)   \int_{\Omega}(a-u_{p,\lambda})^{2\theta-2}|\nabla(a-u_{p,\lambda})|^2\eta_r^{2\tau} dx +c_\delta \int_{\Omega}  (a-u_{p,\lambda})^{2\theta}|\nabla (\eta_r^\tau) |^2dx   ,
\end{eqnarray*}
where $c_\delta>0$ depends on $\delta$.
Then together with (\ref{4.3.2}), we deduce that
\begin{eqnarray*}
 \int_{\Omega}\frac{p\lambda   }{(a-u_{p,\lambda})^{p+1-2\theta}}\eta_r^{2\tau}  dx
  &\le&   \frac{ \theta^2}{1-2\theta}\frac{\theta^2+\delta}{\theta^2-\delta}\int_{\Omega}[ \Delta a+\frac{\lambda}{(a-u_{p,\lambda})^p}] (a-u_{p,\lambda})^{2\theta-1}  \eta_r^{2\tau}dx
   \\[2mm]&&+  \bar c_\delta  \int_{\Omega}(a-u_{p,\lambda})^{2\theta} |\nabla (\eta_r^\tau) |^2dx,
\end{eqnarray*}
where $\bar c_\delta=1+\frac1\delta+ c_\delta$.
Thus,
\begin{eqnarray*}
&&\lambda\Big(p-\frac{\theta^2}{1-2\theta}\frac{\theta^2+\delta}{\theta^2-\delta}\Big)\int_{\Omega}\frac{\eta_r^{2\tau}}{(a-u_{p,\lambda})^{p+1-2\theta}}dx
\\[2mm] &\le&  \frac{\theta^2}{1-2\theta}\frac{\theta^2+\delta}{\theta^2-\delta}  \int_{\Omega} |\Delta a|  (a-u_{p,\lambda})^{2\theta-1} \eta_r^{2\tau}   dx
 +  \bar c_\delta  \int_{\Omega}(a-u_{p,\lambda})^{2\theta} |\nabla (\eta_r^\tau)|^2dx
    \\[2mm] &\le&  C_\delta\Big( \|\Delta a\|_{L^\infty(\Omega_r)}\int_{\Omega} (a-u_{p,\lambda})^{2\theta-1} \eta_r^{2\tau}   dx+\tau^2\|\nabla \eta_r\|_{L^\infty(\Omega_r)} \int_{\Omega}(a-u_{p,\lambda})^{2\theta}\eta_r^{2(\tau-1)}dx\Big).
    \end{eqnarray*}
 
Since $a\in C_0^2(\Omega)$
$$  \|\Delta a\|_{L^\infty(\Omega_r)} \le C_r $$
 for some $C_r>0$ depends on $r$, then letting 
 $$\theta\in(-p-\sqrt{p+p^2},0),$$ we have that $p-\frac{\theta^2}{1-2\theta}>0$ and choose $\delta>0$ such that 
 $$p-\frac{\theta^2}{1-2\theta}\frac{\theta^2+\delta}{\theta^2-\delta}>0,$$
 then by  H\"{o}lder inequality, we obtain that
\begin{eqnarray*}
&& \lambda\Big(p-\frac{\theta^2}{1-2\theta}\frac{\theta^2+\delta}{\theta^2-\delta}\Big)\int_{\Omega}\frac{\eta_r^{2\tau}}{(a-u_{p,\lambda})^{p+1-2\theta}}dx
\\[2mm] & \le & C\Big(  \int_{\Omega} (a-u_{p,\lambda})^{2\theta-1} \eta_r^{2\tau}   dx+ \int_{\Omega}(a-u_{p,\lambda})^{2\theta}\eta_r^{2(\tau-1)}\Big)
\\[2mm] & \le &  C |\Omega|^{ \frac{p}{p+1-2\theta}}\left(\int_{\Omega}\frac{\eta_r^{2\tau}}{(a-u_{p,\lambda})^{p+1-2\theta}}dx\right)^{\frac{1-2\theta}{p+1-2\theta}}
  +
C |\Omega|^{ \frac{p+1}{p+1-2\theta}} \left(\int_{\Omega}\frac{\eta_r^{2\tau}}{(a-u_{p,\lambda})^{p+1-2\theta}}dx\right)^{\frac{-2\theta}{p+1-2\theta}},
\end{eqnarray*}
where $\tau=\frac{p+1-2\theta}{p+1}$ and  
$$2\tau \frac{ -2\theta}{p+1-2\theta}=2(\tau-1).$$
Thus,
there exists $c_{24}>0$ independent on $ \lambda$ such that
$$
 \int_{\Omega}\frac{\eta_r^2}{(a-u_{p,\lambda})^{p+1-2\theta}}dx\le c_{24},
$$
which implies that 
 \begin{equation}\label{4.3.4-0-0}
 \int_{\Omega_{2r}}\Big(\frac1{(a-u_{p,\lambda})^p}\Big)^{q}dx\le c_{24}
 \end{equation}
 for any $r>0$ and $1\leq q<3+\frac1p+2   \sqrt{1+\frac1p}$.
Therefore, passing to the limit as $\lambda\to \lambda^*_{\gamma,p}$, for above $q$ we have that  
 $u_{p,\lambda^*_{\gamma,p}}\in W^{2,q}_{loc} (\Omega)$. 
  \hfill$\Box$\medskip

 \noindent {\bf Proof of Theorem \ref{teo 6}. }  
 $(i)$  Recall that $p^\#_{_N} >0$ be defined in (\ref{ex p}). When $N\geq6$, it is the unique zero point in $(0,+\infty)$ of  
$f_0(p)-N=0$,
where we call $f_0(t) =3+\frac1t  + 2\sqrt{1 +\frac1t}.$

 When $f_0(p) > N$, there exists $\sigma\in[1,2)$ such that 
$$p+1-2\theta\geq q_\sigma$$
 for some $\theta\in\big(-p-\sqrt{p+p^2},0\big)$.
Choosing $\sigma=1$, for $p<p^\#_{_N}$ there holds
  \begin{equation}\label{4.3.4-00}
  \frac1p\Big(p+1+2(p+\sqrt{p+p^2})\Big)> N,
  \end{equation}
  that is $f_0(p)>N$ for $p<p^\#_{_N}$. 
Then by (\ref{4.3.4-0-0})   there exists $q_1\geq N$
  \begin{equation}\label{4.3.4-001}
 \int_{\Omega_{2r}}\Big(\frac1{(a-u_{p,\lambda^*_{\gamma,p}})^{p}}\Big)^{q_1}dx\le c_{24},
\end{equation}

Thus,  for $p\in\I_\gamma$, $ p<p^\#_{_N}$ and  $\lambda=\lambda^*_{\gamma,p}$
it follows by Lemma \ref{lm 4.2} that for $r\in(0,1)$ 
$$\inf_{ \Omega_{2r}} (a-u_{p,\lambda^*_{\gamma,p}})>0$$
and $u_{p,\lambda^*_{\gamma,p}}$ is  $C^{2}_{loc}(\Omega)$ and it is a classical solution of (\ref{eq 1.1})  under our setting.\smallskip

$(ii)$ From Lemma  \ref{lm 8.1}, we have that 
$u_{p,\lambda^*_{\gamma,p}}\in W^{2,q}_{loc}(\Omega)$  for   any 
$$ q\in\Big[1, f_0(p)\Big).$$

If $p^\#_{_N}<p^*_\gamma$, for   
$$
 p^\#_{_N}\leq p< \min\{p^*_\gamma, q^\#_{_N}\}  \quad  {\rm if}\ \  N\geq 7,
$$
 then $\frac{N}{2}<f_0(p)\leq N$ and
 $$\alpha=2-\frac{N}{f_0(p)}\in(0,1].$$
and we have that
$$\frac1{(a-u)^p}\in L^{\frac{N}{2-\alpha}}(\Omega_{2r})$$
and then $u\in W^{2, \frac{N}{2-\alpha}}(\Omega_{3r})$ and by Sobolev's Theorem, we
deduce that $u\in C^\alpha(\Omega_{4r})$. 
 \hfill $\Box$%\medskip

\bigskip

\bigskip

\medskip
 
\noindent{\small {\bf Acknowledgements.} H. Chen is supported by NSFC,  No. 12071189, by Jiangxi Province Science Fund No. 20212ACB211005. 
Y.Wang is supported by NNSF of China, No:12001252, by the Jiangxi Provincial Natural Science Foundation, No:20202ACBL201001.     F. Zhou is supported
by Science and Technology Commission of Shanghai Municipality (STCSM), Grant No. 18dz2271000
and also supported by NSFC, No. 12071189. }

\bibliographystyle{amsplain}

\end{document}